\documentclass[11pt]{amsart}
\usepackage{amsmath,amsthm, amscd, amssymb, amsfonts}

\newcommand\boxe{\begin{tabular}{|p{0,1cm}|}
\hline \\ \hline \end{tabular}}

\newcommand\boxee{\begin{tabular}{|p{0,3cm}|}
\hline \\ \hline \end{tabular}}

\newcommand\boxtru{\begin{tabular}{|p{0,1cm}}
\\ \hline  \end{tabular}}

\newcommand\boxtruu{\begin{tabular}{|p{0,5cm}}
\\ \hline  \end{tabular}}

\newcommand\boxtr{\begin{tabular}{p{0,1cm}|}
\hline \\  \end{tabular}}

\newcommand{\Y}{{\mathcal Y}}

\newcommand{\ku}{\Bbbk}

\newcommand{\X}{{\mathbb X}}

\newcommand{\G}{{\mathcal G}}

\newcommand{\C}{{\mathcal C}}
\newcommand{\D}{{\mathcal D}}
\newcommand{\Ec}{{\mathcal E}}
\newcommand{\uno}{{\bf 1}}

\newcommand{\B}{{\mathcal B}}
\newcommand{\T}{{\mathcal T}}
\newcommand{\Hc}{{\mathcal H}}
\newcommand{\Vc}{{\mathcal V}}
\newcommand{\Pc}{{\mathcal P}}

\newcommand{\Ss}{{\mathcal S}}

\newcommand{\Aut}{\operatorname{Aut}}

\newcommand\Rep{\operatorname{Rep}}

\newcommand\Hom{\operatorname{Hom}}
\newcommand\Opext{\operatorname{Opext}}
\newcommand\Tot{\operatorname{Tot}}

\newcommand{\fde}{{\triangleright}}
\newcommand{\fiz}{{\triangleleft}}

\newcommand{\prin}{t}
\newcommand{\fin}{b}
\newcommand{\pri}{r}
\newcommand{\fine}{l}
\newcommand\rh{\sim_{H}}
\newcommand\rv{\sim_{V}}
\newcommand\rd{\sim_{D}}

\numberwithin{equation}{section}\theoremstyle{plain}

\newtheorem{theorem}{Theorem}[section]
\newtheorem{lema}[theorem]{Lemma}

\newtheorem{prop}[theorem]{Proposition}

\theoremstyle{definition}
\newtheorem{definition}[theorem]{Definition}
\newtheorem{exa}[theorem]{Example}

\theoremstyle{remark}
\newtheorem{obs}[theorem]{Remark}

\newcommand\id{\operatorname{id}}
\newcommand\idd{\mathbf{id}}
\newcommand\iddh{\mathbf{id}}
\newcommand\iddv{\mathbf{id}}

\def\pf{\begin{proof}}
\def\epf{\end{proof}}

\theoremstyle{remark}

\begin{document}

\renewcommand{\baselinestretch}{1.2}
\renewcommand{\thefootnote}{}
\thispagestyle{empty}
\title{Double categories and quantum groupoids}
\author{ Nicol\'as Andruskiewitsch and Sonia Natale}
\address{Facultad de Matem\'atica, Astronom\'\i a y F\'\i sica
\newline \indent
Universidad Nacional de C\'ordoba
\newline
\indent CIEM -- CONICET
\newline
\indent (5000) Ciudad Universitaria, C\'ordoba, Argentina}
\email{andrus@mate.uncor.edu, \quad \emph{URL:}\/
http://www.mate.uncor.edu/andrus}
\email{Sonia.Natale@dma.ens.fr, natale@mate.uncor.edu 
\newline \indent \emph{URL:}\/ http://www.mate.uncor.edu/natale} 
\dedicatory{To Susan Montgomery, on her 60th. birthday}
\thanks{This work was partially supported by CONICET,
Agencia C\'ordoba Ciencia, ANPCyT    and Secyt (UNC)}
\date{September 29, 2003}
\begin{abstract} We give the construction of a class of weak Hopf algebras (or quantum groupoids)
associated to a matched pair of groupoids and  certain cocycle
data. This generalizes a now well-known construction for Hopf
algebras, first studied by G. I. Kac in the sixties. Our approach
is based on the notion of double groupoids, as introduced by
Ehresmann.
\end{abstract}

\maketitle

\section*{Introduction}

An  \emph{exact factorization} of a group $\Sigma$ is a pair of
subgroups $G$, $F$ such that the multiplication map induces a
bijection $m: F\times G \to \Sigma$. Given an exact factorization
of a group $\Sigma$, there are a right action $\fiz: G\times F \to
G$ and a left action $\fde: G\times F \to F$ defined by $sx = (s
\fde x) (s \fiz x)$, for all $s \in G$, $x \in F$. These actions
satisfy the compatibility conditions
\begin{align}\label{comp1}
s \fde xy & = (s \fde x) ((s \fiz x) \fde y), \\
\label{comp2} st \fiz x & = (s \fiz (t \fde x)) (t \fiz x),
\end{align}
for all $s, t \in G$, $x, y \in F$. It follows that $s \fde 1 = 1$ and
$1 \fiz x = 1$, for all $s \in G$, $x \in F$.
Such a data of groups and compatible actions is called a {\it matched
pair} of groups. Conversely, given a matched pair of groups $F$, $G$, one
can find a group $\Sigma$ together with an exact factorization $\Sigma = F G$.

\medbreak Let $\ku$ be a field.
In the early eighties, Takeuchi achieved a construction which, starting from a matched pair $F$, $G$ of finite groups, gives a (in general not commutative and not cocommutative) Hopf algebra $H : = \ku^G \# \ku F$
\cite{tak0}. 
This Hopf algebra fits into an exact sequence
$$\begin{CD}1 @>>> \ku^G @> \iota>>  \ku^G \# \ku F @>\pi>> \ku F @>>> 1\end{CD};$$
$\ku^G \# \ku F $ is called a bismash product; it is semisimple and cosemisimple if the characteristic of $\ku$ is relatively prime to the order of $\Sigma$.
The same construction was also presented by Majid \cite{maj}.
A more general instance of this construction can be done by adjoining a certain cohomological data  associated to the matched pair: namely a pair of 2-cocycles $\sigma: F\times F \to (\ku^G)^{\times}$ and $\tau: G\times G \to
(\ku^F)^{\times}$ satisfying appropriate compatibility conditions. 
In this way,  all Hopf
algebras $H$ which fit into an exact sequence
$\begin{CD}1 @>>> \ku^G @> \iota>>  H @>\pi>> \ku F @>>> 1\end{CD}$ are obtained.
The compatibility conditions have an elegant description in terms of the total complex associated
to a double complex that combines the group cohomologies of $G$,  $F$ and $\Sigma$.
It turns out that this more general construction, and the cohomology theory behind it, had already been discovered by G. I. Kac \cite{k}.
Explicit computations can be done with the help of the so-called Kac exact sequence \emph{loc. cit.}
The  study of this cohomology theory has been later pursued by Masuoka; see the paper \cite{ma-newdir} and references therein for details on this topic.

Quantum groups appart, this construction gave rise to  one of the
first genuine examples of non-commutative non-cocommutative Hopf algebras. 
More recently, it was shown that the Hopf algebras $\mathbb C^G \# \mathbb C F$ 
are exactly those having a positive basis \cite{lyz2}.

\medbreak In our previous paper \cite{AN}, we discussed 
 braided Hopf algebras $R$ which fit into an exact sequence
$\begin{CD}1 @>>> \ku^G @> \iota>>  R @>\pi>> \ku F @>>> 1\end{CD}$.
 The present paper was inspired by a comment of the referee of \cite{AN},
pointing out a pictorial description of the standard basis of $\ku^G \# \ku F$, 
which gives a more compact form to the constructions. 
It turns out that this pictorial description, also present
in \cite{maj-lib, tak, dvvv}, can be stated in the language of double categories. 
These have been introduced by Ehresmann \cite{ehr}.   A double category can be defined as a 
$\C$-structured category, where $\C$ is the category of small categories and functors; that is, as a category 
object in the category of small categories.

Roughly, a small double category consists of a set of 'boxes'  $\B = \{ A, B, \dots \}$, 
each box having colored edges (and vertices) $$A =
\begin{matrix} \quad t \quad \\ l \,\, \boxe \,\, r \\ \quad b\quad
\end{matrix};$$ boxes can be 'horizontally' and 'vertically' composed, both compositions 
subject to a natural interchange law. The description given in the Appendix of \cite{AN} fits  
exactly into this framework: here the categories of vertical and horizontal compositions 
correspond to the transformation groupoids attached to the actions $\fiz: G\times F \to
G$ and $\fde: G\times F \to F$, respectively. In this example the vertical edges of boxes 
are colored by elements of $F$, the horizontal edges are colored by elements of $G$, and 
it has the particularity that every box is uniquely determined by a pair of adjacent edges. 
Moreover, in this case there is only one coloring for the 'vertices' of boxes.

It is then natural to ask: what are the double categories that give rise to Hopf algebras in this fashion? First, we shall consider 
\emph{double groupoids}-- double categories where both the horizontal and vertical compositions are invertible-- to have antipodes.
Now, because of the positive basis Theorem in \cite{lyz2}, we know that the answer should be the double groupoids coming from matched pairs of groups as above.
Still, we can ask: what are the double groupoids that give rise to \emph{weak Hopf algebras}? 

The notion of weak Hopf algebras or quantum groupoids was recently introduced in \cite{bnsz, bsz} as a non-commutative version of groupoids. A relevant feature is that they give rise to tensor categories. A weak Hopf algebra has an algebra and a coalgebra structure; the comultiplication is multiplicative but it does not preserve the unit.

Given a finite double groupoid, we endow the vector space with basis the set of boxes with the groupoid algebra structure of the vertical  groupoid, and with the groupoid coalgebra structure of the horizontal groupoid. We found a sufficient condition to get a quantum groupoid; this is condition (2) in Proposition \ref{vert-hor}.
It turns out that double groupoids satisfying this condition
are equivalent to the \emph{vacant} double groupoids considered by Mackenzie \cite{mk1}: every box be determined by any pair of adjacent edges.  Also, vacant double groupoids are in bijective 
correspondence with matched pairs of groupoids \cite{mk1}. 

Our main result is that every vacant double groupoid gives rise to a  weak Hopf algebra in the way described above; this weak Hopf algebra is semisimple if $\ku$ has characteristic 0.
The corresponding construction is also done by adjoining compatible 2-cocycle data. These 2-cocycle data is part of a first quadrant  double complex, as in the group case; there is as well a "Kac exact sequence for groupoids". We point out that these weak Hopf algebras are involutory. 

Our construction may alternatively be presented without double groupoids, using just exact factorizations of groupoids. We feel however that the language of double categories is not merely accidental; it also appears  again in recent work of Kerler and Lyubashenko on topological quantum field theory \cite{KL}.

\medbreak
The paper is organized as follows. The first section is devoted to double categories and double groupoids. For the convenience of the reader not used to the language of double categories, we include many details and proofs. In the second section we discuss vacant double groupoids. We describe vacant double groupoids in group-theoretical sense in Theorem \ref{doubvacconn}. The third section contains the construction of semisimple quantum groupoids and a presentation of the Kac exact sequence for groupoids. 

\subsubsection*{Acknowledgements} We thank the referee of \cite{AN} for the key remark in her/his report. The results of this paper were reported by the first author at the University of Clermont-Ferrand, June  2003, and at the XV Coloquio Latinoamericano de 
\'Algebra, July 2003.

The research of the second author was done during a postdoctoral stay at the D\' epartement de math\' ematiques et applications of the \' Ecole Normale Sup\' erieure, Paris; she thanks Marc Rosso for the kind hospitality.

Part of the work of the first author was done during a visit 
to the University of Reims (February - May 2003); 
he is very grateful to J. Alev for his generous invitation.

\section{Double categories and double groupoids}

\subsection{Definition of double categories}

\

Let $\C$ be a category with pullbacks. Recall that a category object in  $\C$
(or a category
internal to $\C$) is a collection $(A, O, s, t, \id, m)$,
where $A$ (``arrows") and $O$ (``objects") are objects in $\C$;
$s, e: A \to O$ (``source" and "end = target", respectively), $\id: O \to A$
(``identities") and  $m: A {\,}_{e}\times_{s} A \to A$
(``composition") are arrows in $\C$; subject to the usual
associativity and identity axioms. Similarly, a groupoid object in
$\C$ is a category object in $\C$ with all ``arrows" invertible,
which amounts to the existence of a map $\Ss: A \to A$ with
suitable properties.

\bigbreak
{\bf Notation.} Along this paper, in the case where $f, g$ are composable arrows in a category,
their composition $m(f, g)$ will be indicated by juxtaposition: $m(f, g) = fg$ (and not $gf$).

\begin{definition} A (small) \emph{double category}  $\T$ consists of the following data:

\medbreak
\begin{itemize}
\item Four non-empty sets: $\B$ (boxes),
$\Hc$ (horizontal edges), $\Vc$ (vertical edges) and $\Pc$
(points);

\medbreak
\item eight boundary functions: $\prin, \fin: \B \to \Hc$; \quad
$\pri, \fine: \B \to \Vc$; \quad $\pri, \fine: \Hc \to \Pc$; \quad $\prin,
\fin: \Vc \to \Pc$;

\medbreak
\item four identity functions: $\id: \Pc \to \Hc$; \quad
$\id: \Pc \to \Vc$; \quad $\idd: \Hc \to \B$; \quad $\idd: \Vc \to \B$;

\medbreak
\item four composition functions, all denoted by
$m$:
$$\B {\,}_{\fin}\times_{\prin} \B \to \B \quad\text{(vertical composition),}
\quad \B {\,}_{\pri}\times_{\fine} \B \to \B \quad\text{ (horizontal composition),}$$
$$\Hc {\,}_{\pri}\times_{\fine} \Hc \to \Hc, \quad \Vc {\,}_{\fin}\times_{\prin} \Vc \to \Vc;$$

\end{itemize}
such that  the following axioms are satisfied.

\bigbreak

{\bf Axiom 0.} $(\B, \Hc, \prin, \fin, \idd, m)$, $(\B, \Vc, \fine,
\pri, \idd, m)$, $(\Hc, \Pc, \fine, \pri, \id, m)$, $(\Vc, \Pc,
\prin, \fin, \id, m)$ are categories.

\bigbreak

{\bf Axiom 1.} Four identities between possible functions from $\B$ to
$\Pc$, namely $$tr = rt, \quad tl = lt, \quad bl = lb, \quad br = rb.$$

This axiom allows to depict graphically $A\in \B$ as a box $$A =
\begin{matrix} \quad t \quad \\ l \,\, \boxe \,\, r \\ \quad b\quad
\end{matrix}$$ where $t(A) = t$, $b(A) = b$, $r(A) = r$, $l(A) = l$,
and the four vertices of the square representing $A$ are $tl(A)$,  $tr(A)$, $bl(A)$, $br(A)$.
Of course, $t,b,r$ and $l$ mean, respectively, `top', `bottom', `right' and `left'.
Most of the remaining axioms will be stated in this pictorial representation.

\bigbreak \emph{Warning}. A box $A\in \B$ is, in general, \emph{not} determined
by its four boundaries $t,b,r,l$.

We shall write
$A \vert B$ if $r(A)=l(B)$ (so that $A$ and $B$ are horizontally composable), and
$\displaystyle\frac{A}{B}$ if $b(A)=t(B)$ (so that $A$ and $B$ are vertically composable).

The notation $AB$ (respectively, $\begin{matrix}A
\\B\end{matrix}$) will indicate the horizontal (respectively, vertical) compositions,
whenever $A$ and $B$ are composable in the appropriate sense.

\bigbreak
{\bf Axiom 2.} \emph{Compatibility of the compositions with the boundaries}.

Let $A = \begin{matrix} \quad t \quad
\\ l \,\, \boxe \,\, r  \\ \quad b\quad
\end{matrix}$ and $B = \begin{matrix} \quad u \quad \\ s \,\, \boxe \,\, m \\ \quad c \quad
\end{matrix}$ in $\B$.
\begin{flalign} \label{ax2-1} & \text{If } \quad A \vert B, \quad  \text{ then }  \quad AB =
\begin{matrix} \quad tu \quad \\ l \,\, \boxe \,\, m \\ \quad bc \quad
\end{matrix}, &
\\ \label{ax2-2} & \text{If }   \qquad \displaystyle\frac{A}{B},  \quad \text{ then }
\quad \begin{matrix}A\\B\end{matrix} =
\begin{matrix} \quad t \quad \\ ls \,\, \boxe \,\, rm \\ \quad c \quad
\end{matrix}.&
\end{flalign}

The notation $\begin{tabular}{p{0,4cm}|p{0,4cm}} $A$ & $B$ \\
\hline $C$ & $D$ \end{tabular}$ means that all possible horizontal
and vertical products are allowed; in view of Axiom 2, this implies that
$\displaystyle \frac{AB}{CD}$, $\begin{tabular}{p{0,4cm}|p{0,4cm}}
$A$ & $B$ \\ $C$ & $D$ \end{tabular}$.

\bigbreak
{\bf Axiom 3.} \emph{Interchange law between horizontal and vertical compositions}.
If $\begin{tabular}{p{0,4cm}|p{0,4cm}} $A$ & $B$ \\ \hline $C$ &
$D$ \end{tabular}$, then
\begin{equation}\label{permut}
\begin{matrix} A B \\ C D  \end{matrix} :=
\begin{matrix} \{A B\} \\ \{C D\} \end{matrix} =
\left\{\begin{matrix} A  \\ C   \end{matrix}\right\}
\left\{\begin{matrix} B \\ D  \end{matrix}\right\}.
\end{equation}

A consequence of this Axiom is that, given $r \times s$ boxes $A_{ij}$
with horizontal and vertical compositions allowed as in the following arrangement
$$\begin{tabular}{p{0,8cm}|p{0,8cm}|p{0,8cm}|p{0,8cm}}
$A_{11}$ & $ A_{12}$ & \dots & $A_{1s}$ \\ \hline
$A_{21}$ & $ A_{22}$ & \dots & $A_{2s}$ \\ \hline
\dots & \dots & \dots & \dots \\ \hline
$A_{r1}$ & $ A_{r2}$ & \dots & $A_{rs}$  \end{tabular}\ ,$$
then the product
$$\begin{tabular}{p{0,8cm} p{0,8cm} p{0,8cm} p{0,8cm}}
$A_{11}$ & $ A_{12}$ & \dots & $A_{1s}$ \\
$A_{21}$ & $ A_{22}$ & \dots & $A_{2s}$ \\
\dots & \dots & \dots & \dots \\ \
$A_{r1}$ & $ A_{r2}$ & \dots & $A_{rs}$  \end{tabular}$$
is well defined and can be computed associating in all possible ways.

\bigbreak
{\bf Axiom 4.} \emph{Horizontal and vertical identities}. The identity functions $\idd: \Hc \to \B$ (vertical identity),
$\idd: \Vc \to \B$ (horizontal identity) satisfy $$ \idd(g) =
\begin{matrix} \quad g \quad \\ \id \, {l(g)} \,\, \boxe \,\, \id \, {r(g)} \\
\quad g \quad
\end{matrix}, \quad g \in \Hc; \qquad
\idd(x) = \begin{matrix} \quad \id \, {t(x)} \quad \\ x \,\, \boxee
\,\, x \\ \quad \id \, {b(x)} \quad
\end{matrix}, \quad x\in \Vc.
$$

\bigbreak
Note that, in principle, there is an ambiguity when using the notation $\id \, P$ for an element $P \in \Pc$; however, this ambiguity disappears in the pictorial representation. When necessary, the notation $\id_\Hc : \Pc \to \Hc$ and  $\id_\Vc : \Pc \to \Vc$ for the corresponding identity maps will be used.

\bigbreak
{\bf Axiom 5.} \emph{Horizontal and vertical identities of the identities of the points}.
If $P\in \Pc$, then $$\idd \, {\id_\Hc P} = \idd \, {\id_\Vc P};$$ this
box will be denoted $\Theta_P$.

\bigbreak
{\bf Axiom 6.} \emph{Compatibility of the identity with composition of arrows}.
If $g, h \in \Vc$, $x, y \in \Hc$ are composable, then
$$\left\{\begin{matrix} \idd \, g \\ \idd \, h  \end{matrix}\right\}
= \idd \, gh, \qquad \left\{ \idd \,  x \; \idd  \, y \right\} = \idd \, xy.$$
\end{definition}

\begin{lema}\label{lista-ax} {\bf \cite{bs}.}
A double category is a category object in the category of small categories.
\end{lema}

\pf Let $\T$ be a category object in the category of small categories.
Thus $\T = (\mathcal A, \mathcal O, \prin, \fin, \id, m)$,
where $\mathcal A$ and $\mathcal O$ are small categories,
$\prin, \fin: \mathcal A \to \mathcal O$, $\id: \mathcal O \to \mathcal A$
and  $m: \mathcal A {\,}_{\fin}\times_{\prin} \mathcal A \to \mathcal A$
are functors subject to associativity and identity axioms.

Write $\mathcal A = (\B, \Vc, l, r, \idd, m)$ and $\mathcal O = (\Hc, \Pc, l, r, \id, m)$.
The functors $t, b, \id$ and $m$ correspond, respectively, to maps
$$t, b: \B \rightrightarrows \Hc, \quad \id: \Hc \to \B, \quad m: \B {\,}_{\fin}\times_{\prin} \B \to \B,$$ and
$$t, b: \Vc \rightrightarrows \Pc, \quad \id: \Pc \to \Vc, \quad m: \Vc {\,}_{\fin}\times_{\prin} \Vc \to \Vc.$$
The associativity and identity constraints relating the functors $t, b, \id$ and $m$, correspond to the
fact that $(\B, \Hc, \prin, \fin, \idd, m)$ and  $(\Vc, \Pc,
\prin, \fin, \id, m)$ are categories. In what follows we shall see that the functoriality of $t, b, \id$
and $m$, corresponds to the axioms 1--6.

The functoriality of $t$ and $b$ amount to the identities
\begin{align}\label{a1} & tl  = lt, \quad tr = rt, \quad bl = lb, \quad br = rb;\\
\label{a2} & t m  =  m (t \times t), \quad b m = m (b \times b); \\
\label{a3} & t \id  = \id t; \quad b \id = \id b.
\end{align}
One sees that \eqref{a1} corresponds to Axiom 1, \eqref{a2} corresponds to \eqref{ax2-1} in Axiom 2,
and \eqref{a3} corresponds to the right hand side identity in Axiom 4.

The functoriality of $m$ amounts to
\begin{align}\label{b1} & m(l {\,}_{\fin}\times_{\prin} l) = lm, \quad m (r {\,}_{\fin}\times_{\prin} r) = r m; \\
\label{b2} &  m (m {\,}_{\fin}\times_{\prin} m)  =  m (m \times m); \\
\label{b3} & m (\id {\,}_{\fin}\times_{\prin} \id)  = \id m.
\end{align}

Among these, \eqref{b1} corresponds to identity \eqref{ax2-2} in Axiom 2, \eqref{b2} corresponds to
Axiom 3, and \eqref{b3} corresponds to the left hand side identity in Axiom 6.

The functoriality of $\id$ amounts to
\begin{align}\label{c1} & \id l = l \id, \quad \id r = r \id; \\
\label{c2} &  \id m =  m (\id \times \id); \\
\label{c3} & \id \id  = \id  \id.
\end{align}
Whence, \eqref{c1} corresponds to the left hand side identity in Axiom 4, \eqref{c2} corresponds
to the right hand side identity in Axiom 6, and \eqref{c3} corresponds to  Axiom 5.
The constructions and arguments are reversible, and this finishes
the proof of the lemma. \epf

It is customary to represent a double category in the form of four related
categories $$\begin{matrix} \B &\rightrightarrows &\Hc
\\\downdownarrows &&\downdownarrows \\ \Vc &\rightrightarrows &\Pc \end{matrix}$$ subject to the above axioms.
So that the vertical arrows
$$\begin{matrix} \B \\ \downdownarrows \\ \Vc \end{matrix}, \quad \begin{matrix} \Hc
\\ \downdownarrows \\ \Pc \end{matrix},$$
correspond to the categories $\mathcal A$ and $\mathcal O$ of 'arrows' and 'objects',
respectively, while the horizontal arrows $$\B \rightrightarrows \Hc, \quad \Vc \rightrightarrows \Pc,$$
correspond to the functors $\prin, \fin: \mathcal A \rightrightarrows \mathcal O$.

\begin{obs}
The \emph{transpose} of a double category $\T$ is the double category
$\T^t$ with the same boxes and points as $\T$ but interchanging
the r\^oles of the horizontal and vertical categories. A box $B\in \B$ is denoted $B^t$
when regarded in $\T^t$.
This remark allows to deduce some ``horizontal" statements from ``vertical"
ones (or vice versa), by passing to the transpose double category.
\end{obs}

\subsection{Examples}

\

A $2$-category is just a double category where all 'vertical arrows' are
identities; {\it i.e.}, where every element of $\mathcal V$ is an identity. In this case
the elements of $\mathcal H$ are the \emph{morphisms} of the $2$-category, while the elements
of $\mathcal B$ are the \emph{cells}.

\medbreak
Let $\C$ be a small category. Then the class of all square diagrams in $\C$, that
is all the diagrams $$\begin{CD} X @>x>> Y \\ @VfVV @VVgV \\ Z @>y>>  W\end{CD}$$
is a double category (without assuming commutativity of the diagram!).

\medbreak More examples arise considering all commutative diagrams, or all diagrams
whose vertical, respectively horizontal, arrows live in a fixed subcategory.
For instance, given a group $\Sigma$ and two subgroups $F$ and $G$, one can consider
$\Sigma$ as a category with only one object, and then the double category whose
vertical, respectively horizontal, arrows live in $F$, respectively in $G$.

\medbreak A particular case of the preceding remark is the following example,
which shows a connection between double categories and some constructions
in Hopf algebra theory.

\begin{exa}\label{ap-tak}

\medbreak
Let $\fiz: G\times F \to G$, $\fde: G\times F \to F$, be a matched pair of finite groups.
The notation explained in \cite{maj-lib} (see also \cite{dvvv}, \cite{tak}, \cite[Appendix]{AN}) 
coincides with the pictorial representation of the double category defined in what follows.

Take as $\Pc$ a set with a single element: $\Pc : = \{ * \}$. Put also
$\B: = G \times F$,  $\Hc: = G$ and $\Vc : = F$. We have a double category
$$\begin{matrix} G \times F &\rightrightarrows & G
\\\downdownarrows &&\downdownarrows \\ F &\rightrightarrows &\Pc, \end{matrix}$$
defined as follows:
\item $F \rightrightarrows \Pc$ and $G \rightrightarrows \Pc$ are the  groupoids
associated to the group structure on $F$ and $G$, respectively.
The groupoids $G \times F \rightrightarrows  G$ and $G \times F \rightrightarrows  F$
are the ones corresponding to the actions $\fiz$ and $\fde$, respectively;
see \cite[Example 1.2.a]{renault}. More precisely, we have:

\begin{itemize} \item $G \times F \rightrightarrows  G$ is the category whose objects
are the elements of  $G$, the arrows are the elements of $G \times F$ and the source
and target maps are defined, respectively, by
$$t : = \fiz: G\times F \to G, \qquad b: = p_1: G\times F \to G.$$
The composition $m: (G\times F) {\,}_{\prin}\times_{\fin} (G\times F) \to G\times F$
and the identity map $\id: G \to G \times F$, are determined, respectively, by
$$(g, x) . (h, y) : = (g, xy),  \qquad \id(g) = (g, 1), $$
for all $g, h \in G$, $x, y \in F$, such that $g \fiz x = h$. The inverse map is defined
as $(g, x)^{-1} = (g \fiz x, x^{-1})$.

\item $G \times F \rightrightarrows  F$ is the category whose objects are the elements of
$F$, the arrows are the elements of $G \times F$ and the structure maps are defined by
$$r : = \fde: G\times F \to F, \qquad l: = p_2: G\times F \to F,$$
and composition $m: (G\times F) {\,}_{\prin}\times_{\fin} (G\times F) \to G\times F$ and
identity $\id: G \to G \times F$, are determined by
$$(g, x) . (h, y) : = (hg, x),  \qquad \id(x) = (1, x), $$
for all $g, h \in G$, $x, y \in F$, such that $g \fde x = y$. The inverse map is given
in this case by $(g, x)^{-1} = (g^{-1}, g \fde x)$. \end{itemize} \end{exa}

\subsection{Basic properties}

\

In this section we shall prove some basic properties of double categories and introduce some
terminology that will be of use in later sections.

For an element $A \in \B$, we shall use the notation $A^h$ (respectively, $A^v$)
for the  horizontal (respectively, vertical) inverse of $A$, provided they exist; these are
defined respectively by the relations
$$AA^h = \iddh \,{l(A)}, \quad A^hA = \iddh \,{r(A)}, \qquad \begin{matrix}A \\A^v\end{matrix}
= \iddv \,{t(A)}, \quad \begin{matrix}A^v \\A\end{matrix} = \iddv \,{b(A)}.$$

\begin{lema}\label{ah-av} Let $\T$ be a double category and let $A \in \B$. Suppose that
$A = \begin{matrix} \quad t \quad
\\ l \,\, \boxe \,\, r  \\ \quad b\quad
\end{matrix}$ is invertible with respect to horizontal composition. Then
$t = t(A), b = b(A) \in \Hc$ are invertible and we have
$$A^h = \begin{matrix} \quad t^{-1} \quad
\\ r \,\, \boxee \,\, l  \\ \quad b^{-1}\quad \end{matrix}. $$
Similarly, if $A$ is invertible with respect to vertical composition, then
$l = l(A), r = r(A) \in \Vc$ are invertible and we have
$$A^v = \begin{matrix} \quad b \quad
\\ l^{-1} \,\, \boxe \,\, r^{-1}  \\ \quad t \quad \end{matrix}. $$  \end{lema}

\pf It follows from \eqref{ax2-1}, \eqref{ax2-2} and Axiom 4. \epf

\begin{obs}\label{id-inv} Suppose that  $g \in \Vc$ is invertible. It follows from
axioms 4, 5 and 6, that $\iddh \,g$ is vertically invertible and $(\iddh \,g)^v = \iddh \,{g^{-1}}$.
Also, if $x \in \Hc$ is invertible, then $\iddv \,x$ is horizontally invertible
and $(\iddv \,x)^h = \iddv \,{x^{-1}}$. \end{obs}

\begin{lema}\label{h-inv} Let $\T$ be a double category and let $X, R \in \B$.

\bigbreak
(i) Suppose that $\displaystyle\dfrac{X}{R}$, and that $X$ and $R$ are horizontally invertible.
Then $\begin{matrix}X \\R\end{matrix}$ is horizontally invertible, $\displaystyle\dfrac{X^h}{R^h}$
and $\begin{matrix}X^{h} \\R^{h}\end{matrix} = \left\{\begin{matrix}X \\R\end{matrix}\right\}^{h}$.

\bigbreak
(ii) Suppose that $X \vert R$, and that $X$ and $R$ are vertically invertible. Then $X R$ is vertically
invertible,  $X^v \vert R^v$ and $X^{v}R^{v} = \{ XR \}^{v}$.
\end{lema}

\pf We prove part (i), part (ii) being entirely similar; it can also be deduced from part (i)
by going to the transpose double category.
It is clear that
$\displaystyle\dfrac{X^h}{R^h}$. On the other hand, using the interchange law,
we compute
$$\left\{\begin{matrix}X^{h} \\R^{h}\end{matrix}\right\} \left\{\begin{matrix}X \\R\end{matrix}\right\}
= \left\{\begin{matrix}X^{h} X \\R^{h} R \end{matrix}\right\} = \left\{\begin{matrix}\iddh \, r(X)
\\ \iddh \,r(R) \end{matrix}\right\} = \iddh\, r\left\{\begin{matrix}X \\R\end{matrix}\right\},$$ by
axioms 6 and 2. A similar computation shows that $\left\{\begin{matrix}X \\R\end{matrix}\right\}
\left\{\begin{matrix}X^{h} \\R^{h}\end{matrix}\right\} = \iddh \,l\left\{\begin{matrix}X \\R\end{matrix}\right\}$.
This proves the lemma.  \epf

\bigbreak
\begin{lema} Let $A \in \B$ such that $A$ is horizontally and vertically invertible. Assume in
addition that $A^h$ is vertically invertible and $A^v$ is horizontally invertible.
Then $(A^h)^v = (A^v)^h$. \end{lema}

We shall use the notation $A^{-1} : = (A^h)^v = (A^v)^h$; thus
$A^{-1} = \begin{matrix} \quad b^{-1} \quad
\\ r^{-1}  \,\, \boxee \,\, l^{-1}   \\ \quad t^{-1}\quad \end{matrix}$.

\pf We have $\begin{tabular}{p{0,95cm}|p{0,5cm}} $(A^v)^h$ & $A^v$ \\ \hline $A^h$ & $A$ \end{tabular}$
and also $\begin{tabular}{p{0,95cm}|p{0,5cm}} $(A^h)^v$ & $A^v$ \\ \hline $A^h$ & $A$ \end{tabular}$.
Using the axioms in Lemma \ref{lista-ax}, we compute
$$\left\{\begin{matrix}(A^v)^h & A^v \\ A^h & A \end{matrix}\right\} =
\begin{matrix}\{(A^v)^h  A^v \} \\ \{ A^h  A \} \end{matrix} =
\begin{matrix} \iddh\, r(A)^{-1} \\ \iddh\, r(A) \end{matrix} = \iddh\, (r(A)^{-1}r(A)) = \Theta_{br(A)}. $$
On the other hand, we have
$$\left\{\begin{matrix}(A^h)^v & A^v \\ A^h & A \end{matrix}\right\} =
\left\{\begin{matrix}(A^h)^v \\ A^h \end{matrix} \right\}  \left\{ \begin{matrix} A^v \\  A \end{matrix}\right\}  =
\iddv b(A)^{-1} \iddv b(A)  = \iddv (b(A)^{-1}b(A)) = \Theta_{rb(A)}.$$
Hence we get $\left\{\begin{matrix}(A^v)^h & A^v \\ A^h & A \end{matrix}\right\}
= \left\{\begin{matrix}(A^h)^v & A^v \\ A^h & A \end{matrix}\right\}$.
The horizontal cancellation of $\left\{\begin{matrix}A^v \\ A \end{matrix}\right\}$,
which is licit after Remark \ref{id-inv}, and the vertical cancellation of $A^h$ imply the desired identity. \epf

\subsection{Double groupoids}

\

A \emph{double groupoid} \cite{ehr, bs} is a double category
such that all the four component categories are groupoids. Lemma \ref{ah-av} implies
that a double category $$\begin{matrix} \B &\rightrightarrows &\Hc
\\\downdownarrows &&\downdownarrows \\ \Vc &\rightrightarrows &\Pc \end{matrix}$$ is a
double groupoid if and only if the component categories $\B \rightrightarrows \Hc$
and $\B \rightrightarrows \Vc$ are groupoids.

\medbreak
The transpose of a double groupoid is a double groupoid.

\bigbreak
The following technical result will be needed later.

\begin{lema}\label{l-unit} Let $\T$ be a double groupoid, and let $A, B, C \in \B$.
The following statements are equivalent:

\medbreak
(i) $A \vert B \vert C$ and $ABC = \iddv \,{t(ABC)}$;

\medbreak
(ii) there exist $U, V \in \B$ such that

$$\begin{tabular}{p{0,4cm}|p{0,4cm}|p{0,4cm}} $A$ & $U$ & \quad \\ \hline \quad & $V$ & $C$ \end{tabular}, \qquad
\begin{matrix}U \\ V \end{matrix} = B, \qquad AU = \iddv \,{t(AU)}, \qquad VC = \iddv \,{t(VC)};$$

\medbreak
(iii) there exist $U', V' \in \B$ such that

$$\begin{tabular}{p{0,4cm}|p{0,4cm}|p{0,4cm}} \quad & $U'$ & C
\\ \hline A & $V'$ & \quad \end{tabular}, \qquad\begin{matrix}U' \\ V' \end{matrix} = B, \qquad AV'
= \iddv \,{t(AV')} , \qquad U'C = \iddv \,{t(U'C)}.$$

\medbreak
Moreover, in (ii) and (iii) the elements $U, V, U', V'$ are uniquely determined by $A, B, C$, and we have

$$\begin{tabular}{p{1,2cm}|p{0,4cm}|p{1,2cm}} $A$ & $U$ & $\iddv \,{t(C)}$ \\ \hline $\iddv \,{b(A)}$ & $V$ &
$C$ \end{tabular}, \qquad \text{respectively} \qquad
\begin{tabular}{p{1,2cm}|p{0,4cm}|p{1,2cm}} $\iddv \,{t(A)}$ & $U'$ & C
\\ \hline A & $V'$ & $\iddv \,{b(C)}$ \end{tabular}.$$
\end{lema}

\pf Let $U, V$ as in (ii). The uniqueness of $U$ and $V$ follows from cancellation
properties in a double groupoid. Since  $VC = \iddv \,{t(VC)}$, $l(V) = \id \, {lt(VC)}$,
and on the other hand $r(\iddv \,{b(A)}) = \id \, {rb(A)}$. Now,
$$lt(VC) = lt(V) = lb(U) = bl(U) = br(A) = rb(A),$$ since $A\vert U$. This shows
that $\iddv \,{b(A)} \vert V$.  Similarly, using that $AU = \iddv \,{t(AU)}$, we get
$r(U) = \id \, {rt(AU)}$, and
$$rt(AU) = rb(AU) = rb(U) = rt(V) = tr(V) = tl(C) = lt(C).$$ Since $l(\iddv \,{t(C)})
= \id \, {lt(C)}$, we get $\iddv \,{t(C)} \vert U$. Thus we have
$\begin{tabular}{p{1,2cm}|p{0,4cm}|p{1,2cm}} $A$ & $U$ & $\iddv \,{t(C)}$
\\ \hline $\iddv \,{b(A)}$ & $V$ & $C$ \end{tabular}$, as claimed.
The corresponding facts for $U'$ and $V'$ in (iii) are similarly established.

\medbreak
We shall show that (i) $\Longleftrightarrow$ (ii). The proof of the
equivalence of (i) and (iii) is similar and left to the reader.

\medbreak
(i) $\Longrightarrow$ (ii). Let $x = t(A) t(B) \in \Hc$.  Note that
$$r(A^h) = l(A) = \id \, {lt(AB)} = \id \, {l(x)} = l(\iddv \,x);$$ so that $A^h \vert \iddv \,x$.
Define $U : = A^h \iddv \,x$ and $V: = \begin{matrix} U^v \\ B \end{matrix}$.
To see that $V$ is well defined we compute
$$b(U^v) = t(U) = t(A^h) x = t(A)^{-1} t(A) t(B) = t(B),$$ and thus
$\displaystyle\dfrac{U^v}{B}$. By definition, we have that $A \vert U$,
$\displaystyle\dfrac{U}{B}$, $AU = \iddv \,{t(AU)}$ and $\begin{matrix} U \\ V \end{matrix} = B$. Also,
$$r(V) = r(U^v) r(B) = r((A^h)^v) r(\iddv \,x) r(B) = l(A)^{-1} \id \, {r(x)}
r(B) = l(C),$$  the last equality since $B \vert C$ and $ABC = \iddv \,{t(ABC)}$.
Thus $V \vert C$.

We now observe that
\begin{align*}r(U) & = \id \, {rt(B)} = \id \, {tr(B)} = \id \, {tl(C)} = \id \, {lt(C)}, \\
l(V) & = l(U^v) l(B) = l(A^{-1}) l(\iddv \,x) l(B) = r(A)^{-1} l(B) = \id \, {br(A)}
= \id \, {rb(A)}; \end{align*}
hence $U \vert \iddv \,{t(C)}$ and $\iddv \,{b(A)} \vert V$.
Thus $\begin{tabular}{p{1,2cm}|p{0,4cm}|p{1,2cm}} $A$ & $U$ & $\iddv \,{t(C)}$
\\ \hline $\iddv \,{b(A)}$ & $V$ & $C$ \end{tabular}$. On the other hand $\iddv \,{t(ABC)}
= ABC = \left\{\begin{matrix} A & U & \iddv \,{t(C)} \\ \iddv \,{b(A)} & V & C \end{matrix}\right\}$.
This implies that $VC = \iddv \,{t(VC)}$, and part (ii) follows.

\medbreak
(ii) $\Longrightarrow$ (i). Since $VC = \iddv \,{t(VC)}$, $l(V) = \id \, {lt(VC)}$. Also, since
$B = \begin{matrix} U \\ V\end{matrix}$, and $A \vert U$, we have $l(B) = l(U) l(V) = l(U) = r(A)$;
therefore $A \vert B$. Similarly, the assumptions $C \vert V$ and $AU = \iddv \,{t(AU)}$ imply that
$r(B) = R(U) r(V) = R(V) = l(C)$; hence $B \vert C$.

Finally, the interchange law gives $ABC = \left\{\begin{matrix} A & U & \iddv \,{t(C)}
\\ \iddv \,{b(A)} & V & C \end{matrix}\right\} = \iddv \,{t(ABC)}$, and (i) follows. \epf

The following lemma is dual to Lemma \ref{l-unit}. Its proof is left to the reader.

\begin{lema}\label{l-counit} Let $\T$ be a double groupoid, and let $A, B, C \in \B$.
The following statements are equivalent:

\medbreak
(i) $\displaystyle\dfrac{\displaystyle\dfrac{A}{B}}{C}$ and $\begin{matrix} A\\ B
\\ C \end{matrix}= \iddh \,{l\left(\begin{matrix} A\\ B \\ C \end{matrix}\right)}$;

\medbreak
(ii) there exist $U, V \in \B$ such that

$$\begin{tabular}{p{0,4cm}|p{0,4cm}} $A$ & \quad \\ \hline $U$ & $V$ \\ \hline \quad  &
$C$  \end{tabular}, \qquad UV = B, \qquad \begin{matrix} A\\ U \end{matrix}
= \iddh \,{l\left(\begin{matrix} A\\ U \end{matrix}\right)} , \qquad \begin{matrix} V\\ C \end{matrix}
= \iddh \,{l\left(\begin{matrix} V\\ C \end{matrix}\right)};$$

\medbreak
(iii) there exist $U', V' \in \B$ such that
$$\begin{tabular}{p{0,4cm}|p{0,4cm}}\quad & A \\ \hline $U'$ & $V'$
\\ \hline $C$  & \quad  \end{tabular}, \qquad U'V' = B, \qquad
\begin{matrix} A\\ V' \end{matrix} = \iddh \,{l\left(\begin{matrix} A\\ V' \end{matrix}\right)} , \qquad
\begin{matrix} U'\\ C \end{matrix} = \iddh \,{l\left(\begin{matrix} U'\\ C \end{matrix}\right)}.$$

\medbreak
The elements $U, V, U', V'$ in (ii) and (iii) are uniquely
determined by $A, B, C$, and we have
$$\begin{tabular}{p{1,2cm}|p{1,2cm}} $A$ &  $\iddh \,{r(A)}$
\\ \hline $U$ & $V$ \\ \hline  $\iddh \,{l(C)}$ & $C$  \end{tabular}, \quad \text{respectively} \quad
\begin{tabular}{p{1,25cm}|p{1,25cm}} $\iddh \,{l(A)}$ & $A$
\\ \hline $U'$ & $V'$ \\ \hline  $C$ &  $\iddh \,{r(C)}$ \end{tabular}.$$
\qed \end{lema}

The following result is needed in the proof of Theorem \ref{bicross}.

\begin{lema}\label{l-atp3} The following properties hold in a double groupoid $\T$.

(a)  Let $A, X, Y, Z \in \B$ such that
\begin{equation} \label{tricot1}
\begin{tabular}{p{0,4cm}|p{0,8cm}|p{0,4cm}}\quad & $X^{-1}$
& \quad \\  \hline $X$ & $Y$ & $Z$ \\ \hline \quad  & $Z^{-1}$ & \quad  \end{tabular}\quad .
\end{equation}

Then the following conditions are equivalent:

\begin{align}
\label{tricot2}
XYZ &= A.
\\ \label{tricot25} \left\{ \begin{matrix} X^{-1} \\ Y \\ Z^{-1} \end{matrix} \right\} &= A^{-1}.
\end{align}

(b) The collection $X = A = Z$, $Y = A^h$ satisfies \eqref{tricot1}, \eqref{tricot2} and \eqref{tricot25}.
\end{lema}

\pf Part (b) being straightforward, we prove (a). Suppose that \eqref{tricot2} holds.
By assumption we have $b(X) b(Y) b(Z) = b(A)$, and $b(Y) = t(Z^{-1}) = b(Z)^{-1}$.
Therefore $b(X) = b(A) = t(A^v)$. Similarly, $t(Z) = t(A) = b(A^v)$.
Also, $r(X^{-1}) = l(X)^{-1} = l(A)^{-1} = l(A^v)$, and
$r(A^v) = r(A)^{-1} = r(Z)^{-1} = l(Z^{-1})$. This implies that
\begin{equation} \label{tricot3}
\begin{tabular}{p{0,7cm}|p{0,8cm}|p{0,7cm}}$X^v$ & $X^{-1}$ & $A^v$
\\ \hline $X$ & $Y$ & $Z$ \\ \hline $A^v$  & $Z^{-1}$ & $Z^v$  \end{tabular}.
\end{equation}
We compute in two different ways, using the interchange law:
\begin{align*}
\left\{ \begin{matrix} X^v & X^{-1} & A^v \\ X & Y & Z \\ A^v & Z^{-1} & Z^v \end{matrix} \right\}
&= \begin{matrix} \left\{X^v  X^{-1}  A^v\right\} \\ \left\{X  Y  Z\right\}
\\ \left\{A^v  Z^{-1}  Z^v \right\}\end{matrix}  =
\left\{ \begin{matrix}  A^v \\ A \\ A^v \end{matrix} \right\} = A^v
\\ &= \left\{ \begin{matrix} X^v \\ X \\ A^v \end{matrix} \right\}
\left\{ \begin{matrix}  X^{-1}  \\ Y \\ Z^{-1}  \end{matrix} \right\}
\left\{ \begin{matrix}  A^v \\ Z \\ Z^v \end{matrix} \right\}
= A^v \left\{ \begin{matrix}  X^{-1}  \\ Y \\ Z^{-1}  \end{matrix} \right\} A^v;
\end{align*}
thus $\left\{ \begin{matrix}  X^{-1}  \\ Y \\ Z^{-1}  \end{matrix} \right\} = (A^v)^h = A^{-1}$, as claimed.

Conversely, suppose that \eqref{tricot25} holds. Then, in the transpose category  $\T^t$,
$\begin{tabular}{p{1,2cm}|p{0,6cm}|p{1,2cm}}\quad & $X^t$
& \quad \\  \hline $(X^t)^{-1}$ & $Y^t$ & $(Z^t)^{-1}$ \\ \hline \quad  & $Z^t$ & \quad  \end{tabular}$ \quad
and $(X^t)^{-1}Y^t(Z^t)^{-1} = (A^t)^{-1}$;
thus $\left\{ \begin{matrix} X^t \\ Y^t \\ Z^t \end{matrix} \right\} = A^t$
in $\T^t$ by the preceding; that is, \eqref{tricot2} holds in $\T$.
\epf

\section{Vacant double groupoids}

\subsection{Definition and basic properties}

\

The notion of vacant double groupoids appears in \cite[Definition 2.11]{mk1}.

\begin{definition}\label{def-vac} Let $\T$ be a double groupoid.  We shall say that
$\T$ is \emph{vacant} if for any $g\in \Vc$, $x\in \Hc$
such that $r(x) = t(g)$, there is exactly one  $X \in \B$ such that
$X = \begin{matrix} \quad x \quad \\  \,\, \boxe \,\, g \\ \quad \quad
\end{matrix}$. \end{definition}

We give an alternative description of vacant double groupoids that we have found in the 
course of our research; see condition 2 below. This will be useful in Section \ref{wha}.

\begin{prop}\label{vert-hor} Let $\T$ be a double groupoid. The following are equivalent.

\begin{enumerate}
\item $\T$ is vacant.
\item For all $R, S, P \in \B$ such
that $\displaystyle\dfrac{R}{S}$ and $P \vert\left\{\begin{matrix}R \\S\end{matrix}\right\}$,
there exist unique $X, Y \in \B$ such that $\begin{tabular}{p{0,4cm}|p{0,4cm}} $X$ & $R$ \\ \hline $Y$ &
$S$ \end{tabular}$ and $P = \begin{matrix}X \\Y\end{matrix}$.
\item For any $f\in \Vc$, $y\in \Hc$
such that $l(y) = b(f)$, there is exactly one  $Z \in \B$ such that
$Z = \begin{matrix} \quad  \quad \\  f\,\, \boxe \,\,  \\ \quad  y \quad
\end{matrix}$.
\item For all $T, U, Q \in \B$ such
that $T\vert U$ and $\displaystyle\dfrac{Q}{TU}$,
there exist unique $V, Z \in \B$ such that $\begin{tabular}{p{0,4cm}|p{0,4cm}} $V$ & $Z$ \\ \hline $T$ &
$U$ \end{tabular}$ and $Q = VZ$.
\item For any $f\in \Vc$, $x\in \Hc$
such that $l(x) = t(f)$, there is exactly one  $Z \in \B$ such that
$Z = \begin{matrix} \quad x \quad \\  f\,\, \boxe \,\,  \\ \quad   \quad
\end{matrix}$.
\item For any $g\in \Vc$, $y\in \Hc$
such that $r(y) = b(g)$, there is exactly one  $Z \in \B$ such that
$Z = \begin{matrix} \quad  \quad \\  \,\, \boxe \,\, g \\ \quad  y \quad
\end{matrix}$.
\end{enumerate}
\end{prop}

Note that condition (4) says that the transpose double groupoid $\T^t$ is vacant.

\pf (1) $\implies$ (2). Let $R, S, P \in \B$ such
that $\displaystyle\dfrac{R}{S}$ and $P \vert\left\{\begin{matrix}R \\S\end{matrix}\right\}$.
Let $x = t(P)$, $g= l(R)$ and let $X$ be the unique box of the form
$\begin{matrix} \quad x \quad \\  \,\, \boxe \,\, g  \\ \quad \quad\end{matrix}$.
Set $Y = \begin{matrix} X^v \\ P\end{matrix}$; then clearly
$\begin{tabular}{p{0,4cm}|p{0,4cm}} $X$ & $R$ \\ \hline $Y$ &
$S$ \end{tabular}$ and $P = \begin{matrix}X \\Y\end{matrix}$.
Furthermore, if $X'$, $Y'$ are boxes with these properties then clearly
$X'$ should be of the form $\begin{matrix} \quad x \quad \\  \,\, \boxe \,\, g  \\ \quad \quad\end{matrix}$,
hence  $X' = X$, by the uniqueness condition in (1); \emph{a fortiori} $Y' = Y$.

(2) $\implies$ (1). Let $g\in \Vc$, $x\in \Hc$
such that $r(x) = t(g)$.
Put  $P = \iddv \,x$, $R = \iddh \,g$,
$S = R^v = \iddh \,{g^{-1}}$; by  part (2), there exist
unique $X, Y \in \B$ such that $\begin{tabular}{p{0,4cm}|p{0,4cm}} $X$ & $R$ \\ \hline $Y$ &
$S$ \end{tabular}$ and $P = \begin{matrix}X \\Y\end{matrix}$. Clearly,
$X = \begin{matrix} \quad x \quad \\  \,\, \boxe \,\, g \\ \quad \quad
\end{matrix}$.

Let now $X'$ be of the form
$\begin{matrix} \quad x \quad \\  \,\, \boxe \,\, g  \\ \quad \quad\end{matrix}$.
Let $Y' := \begin{matrix} (X')^v \\ P\end{matrix}$; then $\begin{tabular}{p{0,5cm}|p{0,4cm}} $X'$ & $R$ \\ \hline $Y'$ &
$S$ \end{tabular}$ and $P = \begin{matrix}X' \\ Y' \end{matrix}$.
By the uniqueness in part (2), $X = X'$ (and $Y = Y'$).

(3) $\iff$ (4). This follows from the equivalence (1) $\iff$ (2) for $\T^t$.

(1) $\iff$ (3). If $Z$ is of the form
$\begin{matrix} \quad  \quad \\  f\,\, \boxe \,\,  \\ \quad  y \quad\end{matrix}$,
then $Z^{-1}$ is of the form
$\begin{matrix} \quad y^{-1} \quad \\  \,\, \boxee \,\,f^{-1}  \\ \quad   \quad\end{matrix}$.
This remark implies the desired equivalence. 

The proofs of the equivalences (1) $\iff$ (5) and (1) $\iff$ (6) are analogous,
using $Z^h$ and $Z^v$ respectively. This finishes the proof of the proposition.
\epf

\begin{exa} (i) The double category attached to a matched pair of finite groups
as in in Example \ref{ap-tak} is a vacant double groupoid.

(ii) Let $\G$ be any finite groupoid and consider the
double category $\T$ of commuting square diagrams with vertical arrows in $\G$ but with
horizontal arrows only identities. Then $\T$ is a vacant double groupoid.
\end{exa}

We now complete the information in Lemma \ref{l-atp3}.

\begin{lema}\label{l-atp3-bis} Let $\T$ be a vacant double groupoid and let $A\in \B$.
There exists exactly one collection $X, Y, Z \in \B$ such that
\eqref{tricot1}, \eqref{tricot2} and \eqref{tricot25} hold,
namely $X = Z = A$, $Y = A^h$.
\end{lema}

\pf By Lemma \ref{l-atp3}, the collection  $X = Z = A$, $Y = A^h$ satisfies \eqref{tricot1}, \eqref{tricot2} and \eqref{tricot25}.
On the other hand, suppose that  $X, Y, Z \in \B$ is any collection satisfying
\eqref{tricot1}, \eqref{tricot2} and \eqref{tricot25}. Then $l(X) = l(A)$ 
and $b(X)^{-1} =  t(X^{-1}) = t(A^{-1}) = b(A)^{-1}$, hence $X = A$.
Similarly, $Z = A$ and then necessarily $Y = A^h$.
\epf

We list several technical facts about vacant double groupoids 
that are needed later in this paper. 
The straightforward proof is left to the reader.

\begin{prop}\label{basic-vcnt}
Let $\T$ be a vacant double groupoid, $C\in \B$.

\begin{enumerate}
\item[(i)] If a horizontal (resp. vertical) side of $C$ is an identity, then
$C$ is a vertical (resp. horizontal) identity.

\medbreak
\item[(ii)] The set of pairs of boxes $(A,B)$ such that
$\begin{tabular}{p{0,4cm}|p{0,4cm}} $A$ & $B$ \\ \hline $C$ &  \end{tabular}$, $AB = \id t(AB)$
and $\begin{matrix}A \\C\end{matrix} = \id l\left(\begin{matrix}A \\C\end{matrix}\right)$ is either
$$
\begin{cases} \emptyset, \quad &\text{if $C$ is not a horizontal identity,} \\
\{(\Theta_P, \id x): \;x\in \Hc, \; l(x) = P\}, \quad &\text{if } C = \id g, \; P = t(g).
\end{cases}
$$

\medbreak
\item[(iii)] The set of pairs of boxes $(A,B)$ such that
$\begin{tabular}{p{0,4cm}|p{0,4cm}} & $C$ \\ \hline $A$ & $B$  \end{tabular}$, $AB = \id b(AB)$
and $\begin{matrix} C \\ B\end{matrix} = \id r\left(\begin{matrix} C \\B \end{matrix}\right)$ is either
$$
\begin{cases} \emptyset, \quad &\text{if $C$ is not a horizontal identity,} \\
\{(\id h, \Theta_Q): \;h\in \Vc, \; r(h) = Q\}, \quad &\text{if } C = \id g, \; Q = b(g).
\end{cases}
$$

\medbreak
\item[(iv)] The set of pairs of boxes $(A,B)$ such that
$\begin{tabular}{p{0,7cm}|p{0,4cm}} $A$ & $B$ \\ \hline $B^{-1}$ &   \end{tabular}$, $AB = C$
is either
$$
\begin{cases} \emptyset, \quad &\text{if $C$ is not a horizontal identity,} \\
\{(\; \begin{matrix} \quad  z \\  g\,\, \boxe  \\ \quad \end{matrix}\;,
\begin{matrix} \quad  z^{-1} \quad \\  \,\, \boxe \,\,g  \\ \quad
\end{matrix}): \; z\in \Hc, \; l(z) = t(g)\}, \quad &\text{if } C = \id g.
\end{cases}
$$

\medbreak
\item[(v)] The set of pairs of boxes $(A,B)$ such that
$\begin{tabular}{p{0,4cm}|p{0,7cm}}  & $A^{-1}$ \\ \hline $A$ & $B$   \end{tabular}$, $AB = C$
is either
$$
\begin{cases} \emptyset, \quad &\text{if $C$ is not a horizontal identity,} \\
\{(\; \begin{matrix} \quad  \quad \\  g\,\, \boxe  \\ \quad w^{-1}  \end{matrix}\;,
\begin{matrix} \quad  \quad \\  \,\, \boxe \,\,g  \\ \quad  w \quad
\end{matrix}): \; w\in \Hc, \; r(w) = b(g)\}, \quad &\text{if } C = \id g.
\end{cases}
$$
\end{enumerate}
\qed\end{prop}

\subsection{Matched pairs of groupoids}

\

We shall now give a characterization of vacant double groupoids in terms 
of matched pairs of groupoids. This characterization is due to Mackenzie \cite{mk1}. 

\begin{definition}\label{actiongpds} Let $\G$ be a groupoid with
base $\Pc$ and source and target maps $s, e: \G \rightrightarrows \Pc$. 
Let also $p: \Ec \to \Pc$ be a map. A \emph{left action} of $\G$ on $p$ is a map
$\fde : \G {\,}_e\times_p \Ec \to \Ec$ such that
\begin{flalign}
\label{mp-0}  & p(g\fde x) = s(g),&
\\ \label{mp-1}  & g \fde(h \fde x)  = gh \fde x,&
\\ \label{mp-1.5}  & \iddh \,{p(x)} \, \fde  x  =  x,&
\end{flalign}
for all $g, h \in \G$, $x \in \Ec$  composable in the appropiate sense.

Hence, if $\Ec_b := p^{-1} (b)$, then the action
of $g\in \G$ is an isomorphism $g\fde \underline{\quad}:
\Ec_{t(g)} \to \Ec_{s(g)}$.  This somewhat unpleasant notation 
is originated  by our choice of juxtaposition to denote composition.

Given actions of $\G$ on
$p: \Ec \to \Pc$ and $p': \Ec' \to \Pc$, a map $\phi:\Ec \to \Ec'$
is said to \emph{intertwine} the actions if $p = \phi p'$
and $\phi(g \fde x) = g \fde \phi(x)$, for all $g \in \G$, $x \in \Ec$ such that $e(g) = p(x)$.

An action is \emph{trivial} if there exists a set $X$
such that $\Ec = \Pc\times X$, $p$ is the first projection
and $g\fde (e(g), x) = (s(g), x)$ for all $x\in X$, $g\in \G$.

\medbreak
Similarly, a \emph{right action} of $\G$ on $p$ is a map
$\fiz: \Ec {\,}_p\times_s \G \to \Ec$ such that
\begin{flalign}
\label{mp-0.3}  &p(x\fiz g) = e(g),&
\\ \label{mp-2} & (x \fiz g) \fiz h  = x \fiz gh, &
\\ \label{mp-2.5} & x  \fiz \, \iddv \,{p(x)}  = x, &
\end{flalign}
for all $g, h \in \G$, $x \in \Ec$  composable in the appropiate sense.
\end{definition}

\bigbreak
It is convenient to set the following notation:
a \emph{wide subgroupoid} of a groupoid
$\D$ is a groupoid $\Vc$ provided with a functor $F: \Vc \to \D$ which is the identity
on the objects, and induces inclusions on the hom sets.
In other words, it has the same base, and (perhaps) less arrows.

The next two definitions generalize corresponding notions for finite groups,
\emph{cf.} \cite[Definition 2.14]{mk1}.

\begin{definition}\label{matchpairgpds}  
A \emph{matched pair of groupoids} is a pair of
groupoids $t, b: \Vc \rightrightarrows \Pc$, $l, r: \Hc \rightrightarrows \Pc$,
on the same base $\Pc$, endowed with 
a left action $\fde : \Hc {\,}_r \times_t \Vc \to \Vc$ of $\Hc$ on  
$t: \Vc \to \Pc$, and a right action $\fiz : \Hc {\,}_r \times_t \Vc \to \Hc$ of $\Vc$ on 
$r: \Hc \to \Pc$,  
satisfying
\begin{flalign}
\label{mp-0.7}  &b(x\fde g) = l(x \fiz g),&
\\ \label{mp-3} & x \fde fg  = (x \fde f) ((x \fiz f) \fde g), &
\\ \label{mp-4} & xy \fiz g = (x \fiz (y \fde g)) (y \fiz g), &
\end{flalign}
for all $f, g \in \Vc$, $x, y \in \Hc$ such that the compositions are possible.
\end{definition}

We claim that 
\begin{align}\label{iduno}
x \fde \iddv \,{r(x)} &= \iddv \,{l(x)}, \qquad \text{for all  } x \in \Hc, \\
\label{iddos}
\iddh \,{t(g)} \fiz g &= \iddh \,{b(g)}, \qquad \text{for all  }g \in \Vc.
\end{align}
Indeed, $x \fde \iddv \,{r(x)} = x \fde (\iddv \,{r(x)}\iddv \,{r(x)}) =
(x \fde \iddv \,{r(x)})\, ((x \fiz \, \iddv \,{r(x)}) \, \fde \iddv \,{r(x)}) =
(x \fde \iddv \,{r(x)})\,(x \fde \iddv \,{r(x)})$, by
\eqref{mp-3} and \eqref{mp-2.5}. Since  $t(x \fde \iddv \,{r(x)}) = l(x)$
by \eqref{mp-0},  \eqref{iduno} follows.
Similarly \eqref{iddos} follows from \eqref{mp-4}, \eqref{mp-1.5} and \eqref{mp-0.3}.

\begin{definition}\label{exfactgpds} Let $\D \rightrightarrows \Pc$ be a groupoid. An \emph{exact factorization} of
$\D$ is a pair of wide subgroupoids $\Vc$, $\Hc$, such that for any $\alpha \in \D$, there exist
unique $f\in \Vc$, $y \in \Hc$, such that $\alpha = fy$; that is, if the multiplication map $\Vc {\,}_b\times_l \Hc \to \D$ is a bijection.
\end{definition}

\begin{prop}\label{equiv-matchedpair} \cite[Theorems 2.10 and 2.15]{mk1}
The following notions are equivalent.

\begin{enumerate}
\item Matched pairs of groupoids.
\item Groupoids with an exact factorization.
\item Vacant double groupoids.
\end{enumerate}
\end{prop}
If $\Vc$, $\Hc$ is a matched pair of groupoids, the 
groupoid arising in (2) will be denoted
$\D = \Vc\bowtie \Hc$ and called the \emph{diagonal} groupoid.

\pf (1) $\implies$ (2). Let $\fde : \Hc {\,}_r\times_t \Vc \to \Vc$, $\fiz : \Hc {\,}_r\times_t \Vc \to \Hc$ be a matched pair of groupoids on the same base $\Pc$. 
Let $\D$ be the groupoid on the base $\Pc$, with arrows $\Vc {\,}_b\times_l \Hc$, source and target maps $\alpha, \beta: \Vc {\,}_b\times_l \Hc \rightrightarrows \Pc$ given by
$$\alpha(f, y) = t(f), \qquad \beta(f, y) = r(y),$$
and composition defined by the following rule:
$$(f,y) (h, z) = (f (y\fde h), (y\fiz h) z).$$
We shall denote the arrow corresponding to $(f,y)\in \Vc {\,}_b\times_l \Hc$ by
$\begin{matrix} f\,\, \boxtru \,\,  \\ \quad  y \quad\end{matrix}$.
A straightforward verification shows that $\D$ is indeed a well-defined groupoid.
We identify $\Hc$, resp. $\Vc$, with the arrows of the form
$\begin{matrix}   \iddv \,{l(y)}\,\, \boxtru \,\,  \\ \quad  \quad\quad y \end{matrix}$,
resp.
$\begin{matrix}  f\,\, \boxtruu \,\,  \\ \quad \quad \quad  \iddh \,{b(f)} \quad\end{matrix}$.
Then the pair $\Vc, \Hc$  is an exact factorization of $\D$.

\bigbreak
(2) $\implies$ (1). Let $\D$ be a groupoid and let $\Vc$, $\Hc$ be an exact factorization of
$\D$. Define $\fde : \Hc {\,}_r\times_t \Vc \to \Vc$,   and $\fiz : \Hc {\,}_r\times_t \Vc \to \Hc$,
by the formulas
$$
xg = (x\fde g)(x\fiz g), \qquad (x,g) \in \Hc {\,}_r\times_t \Vc.
$$
The uniqueness of the factorization implies that $\fde$, $\fiz$ are well-defined.
It is not difficult to see that these make $\Vc$, $\Hc$ into a matched pair
of groupoids.

\bigbreak
(3) $\implies$ (1). Let $\T$ be a vacant double groupoid.
Given $g\in \Vc$, $x\in \Hc$
such that $r(x) = t(g)$, we set $x \fiz g := b(X)$, $x \fde g := l(X)$ where
$X \in \B$ is the unique box such that
$X = \begin{matrix} \quad x \quad \\  \,\, \boxe \,\, g \\ \quad \quad
\end{matrix}$. That is,
$X = \begin{matrix} \quad \quad x  \quad \\  x \fde g \, \boxee \,\, g \\ \quad \quad x \fiz g  \quad
\end{matrix}$.
The uniqueness gives at once that $\Vc$, $\Hc$, together  with $\fiz$, $\fde$, is a matched pair of groupoids.

\bigbreak
(1) $\implies$ (3). Let $\Vc$, $\Hc$ be groupoids with the same set of objects $\Pc$, endowed with functions
$\fde$,  $\fiz$. Let $\B := \Hc {\,}_r\times_t \Vc$; we denote
$X = (x,g) \in \Hc {\,}_r\times_t \Vc$ by
$X = \begin{matrix} \quad \quad x  \quad \\  x \fde g \, \boxee \,\, g \\ \quad \quad x \fiz g  \quad
\end{matrix}$. We leave to the reader the verification that this gives rise to a vacant double groupoid.
\epf

\begin{obs} Let $\Vc$, $\Hc$, be an exact factorization of a groupoid $\D$
and let $\T$ be the corresponding vacant double groupoid.
Then $\Hc$, $\Vc$, is also an exact factorization of $\D$; the corresponding vacant double groupoid
is $\T^t$.
\end{obs}

\subsection{Structure of vacant double groupoids}

\

\subsubsection{Structure of finite groupoids}\label{st-ffgpd}
We first briefly recall the well-known structure of finite groupoids.
Let $\G \rightrightarrows \Pc$ be a finite groupoid.
We shall denote by $\G(x,y)$
the set of arrows from $x$ to $y$; the set $\G(x,x)$ will be denoted by  $\G(x)$.

There are two basic examples of groupoids:
\begin{itemize}
\item A finite group $G$, considered as the set of arrows of a category with a single object.
\item An equivalence relation $R$ on a finite set $\Pc$; $s$ and $e$ are respectively the first
and the second projection, and the composition is given by $(x,y)(y,v) = (x,v)$.
We shall denote by $\Pc^2$ the equivalence relation where all the elements of $\Pc$ are related; $\Pc^2$ is called the \emph{coarse} groupoid on $\Pc$.
\end{itemize}
If $\G \rightrightarrows \Pc$ and $\G' \rightrightarrows \Pc'$ are two finite groupoids, then two basic operations are:
\begin{itemize}
\item The disjoint union $\G\coprod \G'$, a groupoid on the base $\Pc \coprod \Pc'$.
\item The direct product $\G\times \G'$, a groupoid on the base $\Pc \times \Pc'$.
\end{itemize}

The structure of any finite groupoid can be described with the help of these
basic examples and operations. Namely, let $\G \rightrightarrows \Pc$ be any finite groupoid and define an equivalence relation on $\Pc$ by
$x\sim y$ iff $\G(x,y) \neq \emptyset$. 
We say that $\G$ is \emph{connected}
if $x\sim y$ for all $x,y\in \Pc$. The opposite case is a 
\emph{group bundle}: this is a groupoid such that $x\sim y$ implies $x=y$.
The trivial group bundle is $\Pc \rightrightarrows \Pc$, $s = e = \id$.

\begin{itemize}
\item If $\G$ is connected, then $\G \simeq \G(x) \times \Pc^2$, where $x$ is any element of $\Pc$.  \end{itemize}

Indeed, choose an arrow $\tau_y \in \G(x,y)$ and define $F: \G \to \G(x) \times \Pc^2$, $F(\alpha) = (\tau_z^{-1}\alpha\tau_y, (y,z))$ if $y = s(\alpha)$, $z = t(\alpha)$. Then $F$ is an isomorphism of groupoids.

\begin{itemize}
\item In general, let $P$ be an equivalence class of $\sim$ and let $\G_P$ be corresponding groupoid  on the base $P$; so that $\G_P(x,y) = \G(x,y)$, for all $x, y \in P$. Then $\G \simeq \coprod_{P \in \Pc/\sim} \G_P$. 
\end{itemize}

These remarks provide the general structure of a finite groupoid.

\bigbreak
\subsubsection{Structure of wide subgroupoids}\label{subgpds}
We now give a description of a wide subgroupoid of a connected groupoid in group-theoretical terms. 
We fix a finite non-empty set $\Pc$, a point $O \in \Pc$, and a finite group $D = : D(O)$. 
Let  $\D = D(O) \times \Pc^2$ be the corresponding connected groupoid.

\begin{lema}\label{subgpds-l} 
There is a  bijective correspondence between the following data:

\begin{enumerate}\item Wide subgroupoids $\Hc$ of $\D$;

\item collections $(\sim_H, (H_P)_{P\in \Pc}, (\overline{d_{PQ}})_{P\sim_HQ})$, where 
\begin{itemize}
\item $\sim_H$ is an equivalence relation on $\Pc$,
\item $H_P$ is a subgroup of $D$, $P\in \Pc$,
\item $\overline{d_{PQ}}$ is an element of $H_P\backslash D / H_Q$ such that for any
representative $d_{PQ}$, the following hold: 
\begin{align}\label{widesubgpd1}
d_{PQ} H_Q &= H_P d_{PQ}, \quad \text{if} \quad P\sim_H Q, \\
\label{widesubgpd2}
d_{PQ}d_{QR} &\in H_P d_{PR}, \quad \text{if} \quad P\sim_H Q \sim_H R, \\
\label{widesubgpd3}
d_{PP} &\in H_P, \quad P\in \Pc.
\end{align}
\end{itemize}
\end{enumerate}
\end{lema}

\pf  For each $P \in \Pc$, fix  $\tau_P \in \D(O, P)$.
The correspondence is not natural; it depends on the choice of the
family $(\tau_P)_{P \in \Pc}$.

(1) $\implies$ (2). The equivalence relation  $\sim_H$ is  defined by
$P \sim_H Q$ iff $\Hc(P, Q) \neq \emptyset$. Given $P\in \Pc$, the subgroup
$H_P$ is defined by $H_P := \tau_P \Hc(P) \tau_P^{-1}$. 
If $P\sim_H Q$, choose $h_{PQ} \in \Hc(P, Q)$ and set 
$d_{PQ} := \tau_P h_{PQ} \tau_Q^{-1}$. Then 
$$
H_Q = \tau_Q \Hc(Q) \tau_Q^{-1} = \tau_Q h_{PQ}^{-1} \Hc(P) h_{PQ} \tau_Q^{-1}
= d_{PQ}^{-1} \tau_P \Hc(P) \tau_P^{-1} d_{PQ} = d_{PQ}^{-1} H_P d_{PQ}.
$$
This proves that condition \eqref{widesubgpd1} is satisfied. Conditions \eqref{widesubgpd2} and \eqref{widesubgpd3} are similarly verified.

If we choose another element $\widetilde{h_{PQ}} \in \Hc(P, Q)$ 
then $\widetilde{d_{PQ}} := \tau_P \widetilde{h_{PQ}} \tau_Q^{-1}$ has the same class in $H_P\backslash D / H_Q$ as $d_{PQ}$.
Clearly, if \eqref{widesubgpd1}, \eqref{widesubgpd2}, \eqref{widesubgpd3}
are true for some choice of representatives of $\overline{d_{PQ}}$ then 
they are true for any choice. This finishes the proof of the first implication.

\medbreak
(2) $\implies$ (1). Define a wide subgroupoid $\Hc$ of $\D$ as follows:
if $P, Q\in \Pc$, then
$$
\D(P, Q) \supseteq \Hc(P,Q) := 
\begin{cases} \emptyset, \qquad\qquad\qquad\, \text{if} \quad P\not\sim_H Q;\\
\tau_P^{-1} H_P d_{PQ} \tau_Q, \quad \text{if} \quad P\sim_H Q.
\end{cases}
$$ 
We have to check that $\Hc$ is stable under composition, inverses and identities.
First, if $P\sim_H Q\sim_H R$ then
\begin{align*}
\Hc(P,Q) \Hc(Q, R) & = (\tau_P^{-1} H_P d_{PQ} \tau_Q)(\tau_Q^{-1} H_Q d_{QR} \tau_R) \\
&= \tau_P^{-1} H_P d_{PQ} d_{QR} \tau_R \\
&= \tau_P^{-1} H_P d_{PR} \tau_R,
\end{align*}
where the first equality is by definition, the second by  \eqref{widesubgpd1} and the third
by \eqref{widesubgpd2}. Next, if $P\sim_H Q$, then
$$
\Hc(P,Q)^{-1} = \tau_Q^{-1}   d_{PQ}^{-1} H_P \tau_P =
 \tau_Q^{-1}H_Q d_{PQ}^{-1}  \tau_P =  \tau_Q^{-1}H_Q d_{QP}  \tau_P
= \Hc(Q, P),
$$
using several times \eqref{widesubgpd1}, \eqref{widesubgpd2} and
\eqref{widesubgpd3}.
Similarly, $\id_P \in \Hc(P,P)$ by \eqref{widesubgpd3}.
The second implication is proved.
\epf


\bigbreak
\subsubsection{Double equivalence relations}\label{der}
Let $\Pc$ be a finite non-empty set.
Let $\rh$, $\rv$ be two equivalence relations on 
$\Pc$. Let $\Vc \rightrightarrows \Pc$ be the groupoid defined 
by the relation $\rv$, let 
$\begin{matrix}\Hc  \\\downdownarrows  \\\Pc\end{matrix}$ 
be the groupoid defined by the relation $\rh$,
and let 
$$ \B = \left\{\begin{pmatrix} P & Q \\ R & S \end{pmatrix}
\in \Pc^{2 \times 2}: P \rh Q, \, P\rv R, \, R\rh S, \, 
Q\rv S\right\}.
$$

 Let $\rd$ be the relation defined as follows:
$P \rd Q$ if there exists $R\in \Pc$ such that $P \rh R$,
$R \rv Q$. We shall sometimes denote this as $\begin{matrix}  {P} \, \frac{\quad}{\quad} \, {R} 
\\  \quad\quad \vert\\
 \quad\quad {Q} \end{matrix}$.

\begin{lema}\label{deqrel} (a). The maps $\B \rightrightarrows \Hc$,
$\begin{matrix}\B \\\downdownarrows  \\\Vc\end{matrix}$ given by
$
\begin{pmatrix} P & Q \\ R & S \end{pmatrix} \rightrightarrows 
\begin{matrix} \begin{pmatrix} P & Q \end{pmatrix}\\ 
\begin{pmatrix}R & S \end{pmatrix} \end{matrix}$, 
$\begin{pmatrix} P & Q \\ R & S \end{pmatrix} \rightrightarrows 
\begin{matrix} \begin{pmatrix} P \\  R \end{pmatrix}, & 
\begin{pmatrix} Q \\ S \end{pmatrix} \end{matrix}$, with eviden composition,
define a double groupoid
$\begin{matrix} \B &\rightrightarrows &\Hc
\\\downdownarrows &&\downdownarrows 
\\ \Vc &\rightrightarrows &\Pc \end{matrix}$;
it will be called a \emph{double equivalence class}.

(b). The relation $\rd$ is an equivalence relation
if and only if it is symmetric.

(c). The double equivalence class is vacant if and only if 
 $\rd$ is an equivalence relation and the following condition holds:
\begin{flalign}\label{equiv-vqc}
&\text{If } R, S \in \Pc, \, R \rh S, \, R \rv S, \text{ then } R=S.&
\end{flalign}

(d). Let $\begin{matrix} \B &\rightrightarrows &\Hc
\\\downdownarrows &&\downdownarrows 
\\ \Vc &\rightrightarrows &\Pc \end{matrix}$ be any \emph{vacant} 
double groupoid
and let $\rh$, $\rv$ be the equivalence relations on 
$\Pc$ defined by $\Hc$, $\Vc$ respectively. Then $\rd$
is an equivalence relation on $\Pc$.
\end{lema}

\pf Part (a) is left to the reader. 

(b). The relation $\rd$ is clearly reflexive.  Assume that  $\rd$ is symmetric. Let $P, Q, T \in \Pc$ such that
$P\rd Q$, $Q\rd T$. Then there exist  $R, S \in \Pc$ such that 
$\begin{matrix}  {P} \, \frac{\quad}{\quad} \, {R} 
\\  \quad\quad \vert\\
 \quad\quad {Q} \end{matrix}$, $\begin{matrix}  {Q} \, \frac{\quad}{\quad} \, {S} 
\\  \quad\quad \vert\\
 \quad\quad {T} \end{matrix}$. But then $\begin{matrix}  {S} \, \frac{\quad}{\quad} \, {Q} 
\\  \quad\quad \vert\\
 \quad\quad {R} \end{matrix}$, \emph{i. e.} $S \rd R$. By symmetry, there exists $V\in \Pc$
such that $\begin{matrix}  {S} \, \frac{\quad}{\quad} \, {Q} 
\\  \vert \, \quad \, \vert\\
{V}  \, \frac{\quad}{\quad} \,  {R} \end{matrix}$. Hence 
$\begin{matrix}  {P} \, \frac{\quad}{\quad} \, {V} 
\\  \quad\quad \vert\\
 \quad\quad {T} \end{matrix}$, \emph{i. e.} $P \rd T$.

(c). Assume that the double equivalence class is vacant. 
If $P\rd Q$,  there exists  $R \in \Pc$ such that 
$\begin{matrix}  {P} \, \frac{\quad}{\quad} \, {R} 
\\  \quad\quad \vert\\  \quad\quad {Q} \end{matrix}$.
By vacancy, there exists $V\in \Pc$ such that $\begin{pmatrix} P & R \\ V & Q \end{pmatrix} \in \B$, that is, we have
$\begin{matrix}  {P} \, \frac{\quad}{\quad} \, {R} 
\\  \vert \, \quad \, \vert\\
{V}  \, \frac{\quad}{\quad} \,  {Q} \end{matrix}$, and thus $Q \rd P$.
It follows that $\rd$ is symmetric and hence an equivalence relation, by (b). Now let
$R, S \in \Pc$ be such that  $\begin{matrix}  {R} \, \frac{\quad}{\quad} \, {S} 
\\  \quad\quad \vert\\  \quad\quad {R} \end{matrix}$.
Then both $\begin{matrix}  {R} \, \frac{\quad}{\quad} \, {S} 
\\  \vert \, \quad \, \vert\\
{R}  \, \frac{\quad}{\quad} \,  {R} \end{matrix}$ and 
$\begin{matrix}  {R} \, \frac{\quad}{\quad} \, {S} 
\\  \vert \, \quad \, \vert\\
{S}  \, \frac{\quad}{\quad} \,  {R} \end{matrix}$
belong to $\B$, and by vacancy $R=S$. Thus \eqref{equiv-vqc} holds.

Conversely, assume that $\rd$ is an equivalence relation and \eqref{equiv-vqc} holds.
Then any $\begin{matrix}  {P} \, \frac{\quad}{\quad} \, {R} 
\\  \quad\quad \vert\\  \quad\quad {Q} \end{matrix}$
can be extended to a box
$\begin{matrix}  {P} \, \frac{\quad}{\quad} \, {R} 
\\  \vert \, \quad \, \vert\\
{V}  \, \frac{\quad}{\quad} \,  {Q} \end{matrix}$ in $\B$, since $\rd$ is an equivalence relation.
If also  $\begin{matrix}  {P} \, \frac{\quad}{\quad} \, {R} 
\\  \vert \, \quad \, \vert\\
{U}  \, \frac{\quad}{\quad} \,  {Q} \end{matrix}$ then clearly
$\begin{matrix}  {U} \, \frac{\quad}{\quad} \, {V} 
\\  \quad\quad \vert\\  \quad\quad {U} \end{matrix}$ and hence
$U=V$ by \eqref{equiv-vqc}.

(d). By vacancy, the relation $\rd$ is symmetric, then apply (b).
\epf

\begin{exa} Let $\Sigma$ be a group,  let $F$, $G$ be subgroups of
$\Sigma$, acting on $\Sigma$ respectively on the left and on the right
by multiplication.
Let $\rh$, $\rv$ be the equivalence relations on 
$\Sigma$ defined by these actions. Then the corresponding 
$\rd$ is an equivalence relation, where the associated partition of $\Sigma$ is that given by the double cosets $FqG$, $q \in \Sigma$. Moreover, 
\eqref{equiv-vqc} holds if and only if $G\cap gFg^{-1} = 1$ for  $g \in \Sigma$. (For instance, if the orders
of $F$ and $G$ are relatively prime).
\end{exa}

\bigbreak
Let $r,s$ be natural numbers. Let $\X_{rs}$
be the double equivalence relation on the set 
$\{1, \dots, r\} \times \{1, \dots, s\}$ and with side relations
$$
(i,j) \sim_H (k,l) \iff i = k, \qquad (i,j) \sim_V (k,l) \iff j = l.
$$

\begin{prop}A vacant, finite,  double equivalence relation
is isomorphic to $\X_{rs}$. \end{prop}

\pf Let $Y_1, \dots, Y_r$ be the classes of $\sim_H$ in
 $\Pc$,
and assume that $Y_1 = \{1, \dots, s\}$. We shall define a bijection
$\phi_i: Y_1 \to Y_i$, $2\le i\le s$ and shall prove that
for $j\in Y_1$ and $k\in Y_i$, $j\sim_V k$ if and only if $k = \phi_i (j)$.
So, let us fix $i$ and set $\phi=\phi_i$. 
Fix $a\in Y_i$. If $j\in Y_1$, there exists a unique $k$ such that
$\begin{matrix} {j} \quad\quad \\
\vert \quad\quad \\ {k} \, \frac{\quad}{\quad} \, {a}\end{matrix}$, by \eqref{equiv-vqc}. Set  $\phi_i (j) = k$.
We claim that $\phi_i$ is bijective.

Indeed, assume that  $\phi_i (j) = k = \phi_i (h)$. Then 
$\begin{matrix}  {j} \, \frac{\quad}{\quad} \, {j} 
\\ \vert \quad\quad \\
{k} \quad\quad\end{matrix}$ and 
$\begin{matrix}  {h} \, \frac{\quad}{\quad} \, {j} 
\\ \vert \quad\quad \\
{k} \quad\quad\end{matrix}$; thus $j = h$ and
$\phi$ is injective. Also, let $k\in Y_i$. Then 
$\begin{matrix} {1} \quad\quad \\
\vert \quad\quad \\ {\phi(1)} \, \frac{\quad}{\quad} \, {k}\end{matrix}$,
hence there exists a unique $j\in Y_1$ such that
$\begin{matrix} {1} &\frac{\quad}{\quad}& \, {j} \\
\vert &\quad & \vert \\ {\phi(1)} \, &\frac{\quad}{\quad}& \, {k}\end{matrix}$;  thus $\phi(j) = k$ and
$\phi$ is surjective. 
\epf

\subsubsection{Structure of vacant double groupoids}

Let $\T_1 = \begin{matrix} \B_1 &\rightrightarrows &\Hc_1
\\\downdownarrows &&\downdownarrows 
\\ \Vc_1 &\rightrightarrows &\Pc_1 \end{matrix}$
and  $\T_2 = \begin{matrix} \B_2 &\rightrightarrows &\Hc_2
\\\downdownarrows &&\downdownarrows 
\\ \Vc_2 &\rightrightarrows &\Pc_2 \end{matrix}$ be  double groupoids.
Then we can define double groupoids
$$
\T_1 \coprod \T_2 = \begin{matrix} \B_1 \coprod  \B_2
& \rightrightarrows &\Hc_1 \coprod \Hc_2
\\\downdownarrows &&\downdownarrows 
\\ \Vc_1 \coprod \Vc_2 &\rightrightarrows &\Pc_1 \coprod \Pc_2 \end{matrix},
\qquad 
\T_1 \times \T_2 = \begin{matrix} \B_1 \times  \B_2
& \rightrightarrows &\Hc_1 \times \Hc_2
\\\downdownarrows &&\downdownarrows 
\\ \Vc_1 \times \Vc_2 &\rightrightarrows &\Pc_1 \times \Pc_2 \end{matrix},
$$
and similarly for families of double groupoids. 
If $\T_i$, $i\in I$, is a family of \emph{vacant} double groupoids
then $\coprod_{i\in I} \T_i$ and  $\times_{i\in I} \T_i$
are also vacant.

\bigbreak
Now, let $\T$ be a vacant double groupoid on the base $\Pc$
$\Pc$, let $\D$ be the corresponding diagonal groupoid
and let $\sim_D$ be the associated equivalence relation.
Then $\T = \coprod_{P\in \Pc / \sim_D} \T_P$ where $\T_P$ 
is the ``full" double groupoid with base the class $P$.

\bigbreak
Let us say that a vacant double groupoid $\T$ is \emph{connected}
if the diagonal relation $\sim_D$ is transitive. The preceding remark
shows that it is enough to consider connected vacant double groupoids.

\medbreak
As in subsection \ref{subgpds} above,
we fix a finite non-empty set $\Pc$, a point $O \in \Pc$, 
and a finite group $D$, and set $\D = D(O) \times \Pc^2$.

\begin{theorem}\label{doubvacconn} 
Let $\Hc$, $\Vc$ be wide subgroupoids of $\D$
associated to data 
$(\sim_H, (H_P)_{P\in \Pc}, (\overline{d_{PQ}})_{P\sim_HQ})$ and 
$(\sim_V, (V_P)_{P\in \Pc}, (\overline{e_{PQ}})_{P\sim_HQ})$, respectively,
 as in Lemma \ref{subgpds-l}.

The following are equivalent:

\begin{enumerate}\item $\D = \Vc \Hc$ is an exact factorization.

\item The  following conditions hold:
\begin{enumerate}
\item[(a)] For all $P, Q\in \Pc$, one has
\begin{equation}\label{doubleclasses}
D = \coprod_{R\in \Pc: P\sim_H R, R\sim _V Q} V_P \, e_{PR} d_{RQ} \, H_Q.
\end{equation}
\item[(b)] For all $P\in \Pc$, $V_P \cap H_P = e$.
\end{enumerate}
\end{enumerate}
\end{theorem}

Note that (a) implies that $\sim_D$ is an equivalence relation on $\Pc$, \emph{cf.} subsection \ref{der} above; this agrees with Lemma 
\ref{deqrel} (d). 

\pf
(1) $\implies$ (2).We show (a). Since $\T$ is vacant, we have
\begin{align*}
D  &= \tau_P \D (P,Q) \tau_Q^{-1} = 
\tau_P \left(\coprod_{R\in \Pc: P\sim_H R, R\sim _V Q} 
\Vc(P,R) \Hc(R,Q) \right) \tau_Q^{-1} 
\\ &=
 \coprod_{R\in \Pc: P\sim_H R, R\sim _V Q} 
 V_P e_{PR} \tau_R \tau_R^{-1} H_R d_{RQ}  =
 \coprod_{R\in \Pc: P\sim_H R, R\sim _V Q} 
 V_P e_{PR}  d_{RQ} H_Q,
\end{align*}
as claimed. We show  (b). Let $g\in V_P \cap H_P$. Then 
$\tau_Pg\tau_P^{-1} \in \Vc(P) \cap \Hc(P)$;
but $\Vc(P) \cap \Hc(P) = \id_P$ since $\T$  is vacant.
Thus $g = e$.

(2) $\implies$ (1). Let $x\in \Hc$, $g\in \Vc$
such that $S := r(x) = t(g)$, and set $P = l(x)$,
$Q = b(g)$. That is, we have 
$\begin{matrix}  x \,  \\ \boxtr \,\, g \end{matrix}$.
Now $\gamma := xg \in \D(P, Q) = \coprod_{R\in \Pc: P\sim_H R, R\sim _V Q} 
\Vc(P,R) \Hc(R,Q)$ by assumption (a).
Thus there exist $R\in \Pc$ (unique!),
$f\in \Vc(P,R)$ and $y\in  \Hc(R,Q)$ such that 
$\gamma = fy$, in other words 
$\begin{matrix} \quad x \quad \\ f \,\, \boxe \,\, g \\ \quad y\quad
\end{matrix} \in \B$. Moreover, assume that also
$\begin{matrix} \quad x \quad \\ h \,\, \boxe \,\, g \\ \quad w\quad
\end{matrix}  \in \B$; note $f\in \Vc(P,R)$ and $y\in  \Hc(R,Q)$. 
Then $z := h^{-1}f = wy^{-1} \in  \Vc(R) \cap  \Hc(R)$;
by hypothesis (b), $z = e$. This implies that $\T = \Vc\Hc$ is an
exact factorization.
\epf

\section{Weak Hopf algebras arising from a vacant double groupoid}\label{wha}

Let $F$, $G$ be a matched pair as in Example \ref{ap-tak}.
As explained in many places, see \emph{e. g.} \cite[5.3]{AN}, a bicrossed product Hopf algebra
$\ku^G {}^{\tau}\#_{\sigma} \ku F$ admits a convenient realization in the vector space
with basis $\B$.
In this section we shall discuss
a generalization of this construction.

\subsection{Weak Hopf algebras (quantum groupoids)}

\

Recall \cite{bnsz, bsz} that a \emph{weak bialgebra} structure on a vector space
$H$ over a field $\ku$
consists of an associative algebra structure $(H, m, 1)$, a coassociative coalgebra
structure $(H, \Delta, \varepsilon)$, such that the following are satisfied:
\begin{align}\label{d-mult}\Delta(ab) &= \Delta(a) \Delta(b), \qquad \forall a, b \in H.
\\ \label{ax-unit} \Delta^{(2)} (1 ) &= \left( \Delta(1) \otimes 1 \right)
\left( 1 \otimes \Delta(1) \right) = \left( 1 \otimes \Delta(1) \right)
\left( \Delta(1) \otimes 1 \right).
\\ \label{ax-counit} \varepsilon(abc) &= \varepsilon(ab_1)\varepsilon(b_2c)
= \varepsilon(ab_2)\varepsilon(b_1c), \qquad \forall a, b, c \in H.
\end{align}

A weak bialgebra $H$ is called a \emph{weak Hopf algebra} or a \emph{quantum groupoid}
if there exists a linear map $\Ss: H \to H$ satisfying
\begin{align}\label{atp-1} m(\id \otimes \Ss) \Delta (h) &
= (\varepsilon \otimes \id) \left(  \Delta(1) (h \otimes 1)\right) =: \varepsilon_t(h), \\
\label{atp-2} m(\Ss \otimes \id) \Delta (h) &
= (\id \otimes \varepsilon) \left( (1 \otimes h) \Delta(1) \right)=: \varepsilon_s(h),\\
\label{atp-3} m^{(2)}(\Ss \otimes \id \otimes \Ss) \Delta^{(2)} & = \Ss,
\end{align}
for all $h \in H$. The maps $\varepsilon_s$, $\varepsilon_t$ are respectively called
the source and target maps; their images are respectively called the
source and target subalgebras.
See \cite{nik-v} for a survey on quantum groupoids.
It is known that a weak Hopf algebra is a true Hopf algebra if and only if
$\Delta(1) = 1 \otimes 1$.

\subsection{Weak Hopf algebras arising from vacant double groupoids}

\

Let $\T$ be a  \emph{finite} double groupoid, that is,  $\B$, $\Vc$, $\Hc$ and $\Pc$ are finite sets.

Let $\ku$ be a field (*)\footnote{(*) Most of the constructions in this section are valid over an
arbitrary commutative ring.} 
and let $\ku\T$ denote the  $\ku$-vector space with basis $\B$ together
with the following structures.

\bigbreak \emph{Algebra structure.} Consider the groupoid algebra structure on
$\ku\T$ corresponding to the groupoid $\B \rightrightarrows \Hc$.
Thus the multiplication in $\ku\T$ is given
by $$A.B = \begin{cases} \begin{matrix}A \\B\end{matrix}, \quad
\text{if } \displaystyle\frac{A}{B}, \\
0 , \quad \text{otherwise}, \end{cases}$$
for all  $A, B \in \B$. This multiplication is associative and has a unit $\uno :
= \sum_{x \in \Hc} \iddv \,{x}$. We shall also consider, for any $P\in \Pc$, the elements
$$
{}_{P}\uno= \sum_{x \in \Hc, l(x) = P} \iddv \,{x}, \qquad \uno_{P}= \sum_{x \in \Hc, r(x) = P} \iddv \,{x}.
$$
Clearly, ${}_{P}\uno\; {}_{Q}\uno = \delta_{P, Q}\;  {}_{P}\uno$,
$\uno_{P}\uno_{Q} = \delta_{P, Q} \uno_{P}$, for all $P, Q\in \Pc$.
Hence the subalgebras $\ku \T_s$, respectively $\ku \T_t$,  generated by
$\uno_{P}$, $P\in \Pc$, respectively by ${}_{P}\uno$, $P\in \Pc$, are commutative of dimension
$\vert \Pc \vert$.

\bigbreak \emph{Coalgebra structure.} Dually, we consider the coalgebra structure
on $\ku \T$ dual to the algebra structure of the groupoid algebra corresponding to
the groupoid $\B \rightrightarrows \Vc$. This means that the comultiplication of
$\ku\T$ is determined by $$\Delta(A) = \sum B\otimes C, \qquad A \in \B, $$
where the sum runs over all $B, C$ with $B \vert C$ and $A = B C$.
This comultiplication is coassociative and has counit
$\varepsilon: \ku\T \to \ku$ given by $$\varepsilon(A) =
\begin{cases} 1, \quad \text{if } A = \iddh \,{l(A)},  \\ 0 , \quad \text{otherwise}.
\end{cases}$$

\begin{theorem}\label{bicross} Suppose that $\T$ is a vacant double groupoid. Then
$\ku\T$ is a quantum groupoid with antipode defined by $\mathcal S (A) = A^{-1}$, $\forall A \in \B$.
The source and target subalgebras are respectively $\ku \T_s$, $\ku \T_t$.\end{theorem}

\pf We first show that \eqref{d-mult} holds, using the vacancy of $\T$.
Let $A, B \in \B$. It follows from the definitions that
$$\Delta(A) . \Delta(B)
= \sum_{} \begin{matrix}X \\R\end{matrix} \otimes \begin{matrix}Y \\S\end{matrix},$$
where the sum runs over all elements $X, Y, R, S \in \B$, such that
$\begin{tabular}{p{0,4cm}|p{0,4cm}} $X$ & $Y$ \\ \hline $R$ & $S$ \end{tabular}$,
$XY = A$ and  $RS = B$.
It is thus clear that $\Delta(A) . \Delta(B) = 0 = \Delta(A . B)$, if $A$ and $B$
are not vertically composable. So assume that $\displaystyle\frac{A}{B}$. Hence we have
$$\Delta(A) . \Delta(B) = \sum_{\frac{X}{R}, X^h|A, R^h|B}
\begin{matrix}X \\R\end{matrix} \otimes \begin{matrix}X^{h}A \\R^{h}B\end{matrix}
= \sum_{\frac{X}{R}, X^h|A, R^h|B} \begin{matrix}X \\R\end{matrix}
\otimes \left\{\begin{matrix}X \\R\end{matrix}\right\}^h \left\{\begin{matrix}A \\B\end{matrix}\right\}
= \sum_{U \in \B} U \otimes U^{h} \left\{\begin{matrix}A \\B\end{matrix}\right\},$$
the second equality by Lemma \ref{h-inv}, and the last equality from the substitution
$U = \begin{matrix}X \\R\end{matrix}$.  Using the compatibility relation \eqref{permut}
and the vacancy of $\T$, one sees that the last expression equals $\Delta(A . B)$;
so $\Delta$ is multiplicative as claimed.

We next prove the relationships \eqref{ax-unit} and \eqref{ax-counit}.

We have $\Delta^{(2)} (\uno ) = \sum_{A|B|C, \, ABC = \iddv \,{t(ABC)}} A \otimes B \otimes C$.
On the other hand, we have
$$\left( \Delta(\uno) \otimes \uno \right) \left( \uno \otimes \Delta(\uno) \right)
= \sum_{} A \otimes \begin{matrix} U \\ V \end{matrix} \otimes C,$$
where the sum runs over all $U, V \in \B$ such that
$\begin{tabular}{p{0,4cm}|p{0,4cm}|p{0,4cm}} $A$ & $U$ & \quad \\ \hline \quad & $V$ & $C$ \end{tabular}$,
$\begin{matrix}U \\ V \end{matrix} = B$, $AU = \iddv \,{t(AU)}$, $VC = \iddv \,{t(VC)}$.
This coincides with $\Delta^{(2)} (\uno )$
thanks to the equivalence (i) $\Longleftrightarrow$ (ii) in Lemma \ref{l-unit};
the equivalence (i) $\Longleftrightarrow$ (iii) in the same lemma gives the
equality $\Delta^{(2)} (\uno ) = \left( \uno \otimes \Delta(\uno) \right)
\left( \Delta(\uno) \otimes \uno \right)$. This establishes \eqref{ax-unit}.
Dually, \eqref{ax-counit} follows from Lemma \ref{l-counit}.

\medbreak
We next consider the axioms of the antipode.
We first treat \eqref{atp-1} and \eqref{atp-2}.
Using Proposition \ref{basic-vcnt} (ii), we see that
$$
\varepsilon_t(C) = \begin{cases} 0 \quad &\text{if $C$ is not a horizontal identity,} \\
{}_P\uno, \quad &\text{if } C = \id g, \; P = t(g),
\end{cases} \quad \text{for all } C\in \B.
$$
This coincides with the left hand side of \eqref{atp-1}, by Proposition \ref{basic-vcnt} (iv).
Similarly, using Proposition \ref{basic-vcnt} (iii), we see that
$$
\varepsilon_s(C) = \begin{cases} 0 \quad &\text{if $C$ is not a horizontal identity,} \\
\uno_P, \quad &\text{if } C = \id g, \; P = b(g),
\end{cases} \quad \text{for all } C\in \B.
$$
This coincides with the left hand side of \eqref{atp-1}, by Proposition \ref{basic-vcnt} (v).
Finally, the relation \eqref{atp-3} is equivalent to the identity
$$\sum_{X, Y, Z} \left\{ \begin{matrix} X^{-1} \\ Y \\ Z^{-1} \end{matrix} \right\}
= A^{-1},$$ for all $A \in \B$, where the sum runs over all $X, Y$
and $Z$ in $\B$ such that
$\begin{tabular}{p{0,4cm}|p{0,8cm}|p{0,4cm}}\quad & $X^{-1}$ & \quad
\\ \hline $X$ & $Y$ & $Z$ \\ \hline \quad  & $Z^{-1}$ & \quad  \end{tabular}$
and $XYZ = A$. By Lemmas \ref{l-atp3} and \ref{l-atp3-bis}, the left hand side equals $A^{-1}$. \epf

\bigbreak
By construction, $\ku\T$ is the groupoid algebra of the vertical groupoid $\B \rightrightarrows \Hc$. Using the description in \ref{st-ffgpd}, this groupoid is isomorphic to $\coprod_{H \in \Hc / \sim_h} \B_H$, where $\sim_h$ is the equivalence relation in $\Hc$ defined by $\B \rightrightarrows \Hc$ and $\B_H \rightrightarrows H$ is the connected groupoid $\B(x) \times H^2$ on the class $H$, $x \in H$. Therefore, there is an isomorphism of algebras
\begin{equation}\label{algs}\ku\T \simeq \oplus_{H \in \Hc / \sim_h} \ku \B(x) \otimes M_{n(H)}(\ku), \end{equation}
where $\ku \B(x)$ is the group algebra, and $n(H) = \vert H \vert$.

Similarly, we have an isomorphism of coalgebras
\begin{equation}\ku\T \simeq \oplus_{V \in \Vc / \sim_v} \ku \B(g) \otimes M_{m(V)}(\ku), \end{equation}
where $\sim_v$ is the equivalence relation defined by the horizontal groupoid 
$\B \rightrightarrows \Vc$, 
$\ku \B(g)$ is the group algebra of $\B(g)$, $g \in V$, and $m(H) = \vert V \vert$.

\begin{exa} Suppose that $\ku \T$ is simple as an algebra. Then

(1) $\B \rightrightarrows \Hc$ is the coarse groupoid on $\Hc$;

(2) $\Hc \rightrightarrows \Pc$ is a trivial group bundle;

(3) $\Vc \rightrightarrows \Pc$ is the coarse groupoid on $\Pc$;

(4) $\B \rightrightarrows \Vc$ is a trivial group bundle.

This means that as a weak Hopf algebra $\ku \T$ is the groupoid algebra of the vertical groupoid $\B \rightrightarrows \Hc$, and in particular it is cocommutative. \end{exa}

\pf By \eqref{algs}, $\ku \T$ is simple iff $\vert \Hc / \sim_h\vert = 1$ and $\B(x)$
is trivial for any $x\in \Hc$. Hence (1). If $x: P \to Q$ is in $\Hc$, then there is a box $B$
connecting it to $\id  P$; that is, there exists $g\in \Vc$ such that 
$B =\begin{matrix} \quad \quad \id  \quad \\  x \fde g \, \boxee \,\, g \\ \quad \quad x  \quad
\end{matrix}$ with $x = \id P \fiz g = \id Q$ by \eqref{iddos}.  Hence (2). 
If $P, Q \in \Pc$, there is a unique box  connecting  $\id  P$ and $\id Q$; the vertical
sides connect $P$ and $Q$, hence (3). Finally, let $g: P\to Q$ and $h:R\to S$
be in $\Vc$ if there is box $\begin{matrix} \quad \quad    \quad \\  h \, \boxee \,\, g \\ \quad \quad   \quad
\end{matrix} \in \B$, then $P = Q$ and $R = S$ by (1); and then $g = h$ since 
$\B(\id)$ is trivial; hence (4). \epf

\bigbreak
The proofs of the following statements are straightforward and are left to the reader.

\begin{prop}  Let $\T_1$, $\T_2$ be finite vacant double groupoids. Then there
are isomorphisms of quantum groupoids 
$\ku (\T_1 \coprod \T_2) \simeq \ku \T_1 \times \ku \T_2$,
$\ku (\T_1 \times \T_2) \simeq \ku \T_1 \otimes \ku \T_2$.
\qed \end{prop}

\begin{prop}  Let $\T$ be a finite vacant double groupoid and assume that $\ku 
= \mathbb C$. Then $\mathbb C \T$ is a $C^*$ quantum groupoid \cite{bnsz},
with the involution uniquely defined by $A^* = A^{-1}$, $A\in \B$.
\qed \end{prop}

\subsection{Extensions with cocycles}

\

We begin by recalling the following definition.

\begin{definition}  Let $s, e: \G \rightrightarrows \Pc$ be a groupoid. A \emph{normalized 2-cocycle}
on $\G$ with values in $\ku^{\times}$ is a function
$\sigma: \G {\,}_{s}\times_{e} \G \to \ku^{\times}$ such that
\begin{flalign}
\label{cociclo-sigma} 
& \sigma(\alpha, \beta) \sigma(\alpha\beta, \gamma) 
= \sigma(\beta, \gamma) \sigma(\alpha,\beta\gamma); & \\
\label{norm-sigma} 
& \sigma(\alpha, \id \, {t(\alpha)}) = \sigma(\id \, {b(\alpha)}, \alpha) = 1,&
\end{flalign} 
for all composable $\alpha,\beta,\gamma \in \G$.
\end{definition}

This definition fits into the general framework of groupoid cohomology, due
to Westman \cite{we},  see \cite{renault}.

\bigbreak
Let now $\T$ be a double groupoid. A \emph{normalized  vertical 2-cocycle} is
a 2-cocycle on the groupoid $\B \rightrightarrows \Hc$;
similarly, a \emph{normalized  horizontal 2-cocycle} is
a 2-cocycle on the groupoid $\B \rightrightarrows \Vc$. Thus,
a normalized  vertical 2-cocycle is a function $\sigma$ on the set
of all pairs $(A, B)$ with $\displaystyle \frac{A}{B}$ with values in $\ku^{\times}$
such that
\begin{flalign} \label{cociclo-sigma-bis} & \text{If } \displaystyle
\frac{\displaystyle\frac{A}{B}}{C}, \text{ then } \sigma(A, B)
\sigma(\begin{matrix}A \\B\end{matrix},C) = \sigma(B, C) \sigma(A,
\begin{matrix}B \\C\end{matrix}). & \\
\label{norm-sigma-bis} & \text{If } A \text{ or } B \text{ is a
vertical identity, then } \sigma(A, B) = 1.&
\end{flalign}
Letting $A = B^v= C$, we deduce that
\begin{flalign}
\label{cond-sigma-bis} & \sigma(A, A^v)  =\sigma(A^v, A).& 
\end{flalign}

Analogously, a normalized  horizontal 2-cocycle is a function
$\tau$ on the set of all pairs $(A, B)$ with $A \vert B$,
such that
\begin{flalign}
\label{cociclo-tau-bis} & \text{If } A\vert B\vert C, \text{ then
} \tau(A, B) \tau(AB, C) =\tau(B, C) \tau(A, BC).& \\
\label{norm-tau-bis} & \text{If } A \text{ or } B \text{ is a
horizontal identity, then } \tau(A, B) = 1.&
\end{flalign}
Letting $A = B^h = C$, we deduce that
\begin{flalign}
\label{cond-tau-bis} & \tau(A, A^h)  =\tau(A^h, A).& 
\end{flalign}

\begin{definition} \label{double-cocycle} Let $\T$ be a double groupoid. A \emph{normalized 2-cocycle}
on $\T$ with values in $\ku^{\times}$ is a pair $(\sigma, \tau)$, where $\sigma$ is
a normalized  vertical 2-cocycle, $\tau$ is normalized  horizontal 2-cocycle, and
the following property holds:
\begin{equation}
\label{cociclo-sigma-tau}   \text{If }\quad \begin{tabular}{p{0,4cm}|p{0,4cm}} A & B
\\ \hline C & D \end{tabular},  \quad\text{ then  }\quad \sigma(AB,
CD) \tau \left(\begin{matrix} A \\ C \end{matrix},
\begin{matrix}B \\ D \end{matrix}\right)
= \tau(A, B) \tau(C, D) \sigma(A, C) \sigma(B, D).  \end{equation}
\end{definition}

\bigbreak
Given a normalized vertical 2-cocycle $\sigma$ and a normalized horizontal 2-cocycle $\tau$ on the double groupoid $\T$, one may consider the $\sigma$-\emph{twisted} groupoid algebra structure and, dually, the $\tau$-\emph{twisted} groupoid coalgebra structure on the vector space $\ku \T$ with basis $\B$. The following theorem asserts that, provided that $\T$ is vacant, the compatibility condition \eqref{cociclo-sigma-tau} guarantees that these two structures combine into a weak Hopf algebra structure. 

\bigbreak
\begin{theorem}\label{concociclos} 
Let $\T$ be a vacant double groupoid and let $(\sigma, \tau)$
be a normalized 2-cocycle on $\T$ with values in $\ku^{\times}$. 

(i) Let $\ku{}^{\tau}_{\sigma} \T$ be the vector space with basis $\B$ and 
multiplication and  comultiplication defined, respectively  by
\begin{itemize}\item $A.B = \sigma(A, B) \begin{matrix}A
\\B\end{matrix}$,  if $\displaystyle \frac{A}{B}$,  and $0$ otherwise.
\item $\Delta(A) = \sum \tau (B, C) B\otimes C$ , where the sum is
over all pairs $(B, C)$ with $B\vert C$ and $A = BC$.
\end{itemize}
Then $\ku{}^{\tau}_{\sigma} \T$ is a quantum groupoid with antipode defined by
\begin{equation}\label{antipoda}
\mathcal S (A) = \tau(A, A^h)^{-1} \, \sigma(A^{-1}, A^h)^{-1}\, A^{-1}. 
\end{equation}

The source and target subalgebras are, respectively, the subspaces spanned by $(\uno_P)_{P \in \Pc}$ and $({}_P\uno)_{P \in \Pc}$; so they are commutative of dimension $\vert \Pc \vert$.

(ii) Let $(\nu, \eta)$ be another normalized 2-cocycle on $\T$ with values in $\ku^{\times}$.
Let $\psi: \T \to \ku^{\times}$ be a map and let $\Psi: \ku{}^{\tau}_{\sigma} \T \to \ku{}^{\eta}_{\nu} \T$
be the linear map given by $\Psi(B) = \psi (B) B$, $B\in \B$.
Then $\Psi$ is an isomorphism of quantum groupoids if and only if
\begin{align} \psi\left(\begin{matrix}A
\\B\end{matrix}\right) \sigma(A, B) &= \psi(A)\psi(B) \nu(A, B),
\quad \text{for all } A, B  \in \B \text{ such that  }\frac{A}{B};
\\
\psi(CD) \eta(C, D) &= \psi(C)\psi(D) \tau(C, D),\quad 
\text{for all }  C, D \in \B \text{ such that  } C\vert D.\end{align}
 \end{theorem}

\pf (i). Straightforward computations show that the multiplication
is associative with unit 1 (because of the cocycle and
unitary conditions on $\sigma$), and that the comultiplication
is coassociative with counit $\varepsilon$ (because of the cocycle and
unitary conditions on $\tau$). By \eqref{cociclo-sigma-tau}, $\Delta$ is multiplicative.

We now prove that $\Delta^{(2)} (1 ) = \left( \Delta(1) \otimes 1 \right)
\left( 1 \otimes \Delta(1) \right)$. Arguing as in the proof of Theorem \ref{bicross},
it is enough to check that
$$
\tau \left(A, \begin{matrix}U \\ V \end{matrix}\right)
\tau \left(A \left\{\begin{matrix}U \\ V \end{matrix} \right\}, C\right) =
\tau(A, U)\sigma(U,V) \tau(V, C),$$
if $\begin{tabular}{p{1,2cm}|p{0,4cm}|p{1,2cm}} $A$ & $U$ & $\iddv \,{t(C)}$ \\ \hline $\iddv \,{b(A)}$ & $V$ &
$C$ \end{tabular}$, $AU = \iddv \,{t(AU)}$, $VC = \iddv \,{t(VC)}$.

Now 
 $\begin{tabular}{p{0,4cm}|p{0,4cm}} $A$ & $U$  \\ \hline  $\id$ & $V$ \end{tabular}$
implies $\tau\left(A, \begin{matrix}U \\ V \end{matrix} \right) = \tau(A,U) \sigma(U,V)$,
and $A \left\{\begin{matrix}U \\ V \end{matrix} \right\}= \begin{matrix} AU \\ V \end{matrix} = V$,
and the claimed identity holds.

The proof of $\Delta^{(2)} (1 ) = \left( 1 \otimes \Delta(1) \right)\left( \Delta(1) \otimes 1 \right)$
is similar, using Lemma \ref{l-unit} (iii); so is the proof of \eqref{ax-counit}
using Lemma \ref{l-counit}.

\medbreak
The proof of \eqref{atp-1} and \eqref{atp-2} is as in the proof of Theorem \ref{bicross};
\eqref{cond-sigma-bis} and \eqref{cond-tau-bis} are needed. Note that the source and target maps coincide with those of $\ku \T$.

\medbreak
We now prove \eqref{atp-3}. Given $A\in \B$, we compute:
\begin{align*}
m^{(2)}(\Ss \otimes \id \otimes \Ss) \Delta^{(2)} (A) &= m^{(2)}(\Ss \otimes \id \otimes \Ss) 
\left(\sum_{X\vert Y \vert Z, XYZ = A} \tau(X, Y)\tau(XY, Z) X \otimes Y \otimes Z \right) \\
&= \sum \tau(X, Y)\tau(XY, Z) \tau(X, X^h)^{-1} \, \sigma(X^{-1}, X^h)^{-1}
\tau(Z, Z^h)^{-1} \, \sigma(Z^{-1}, Z^h)^{-1} \\ &\qquad \qquad \times \sigma(X^{-1}, Y) 
\sigma(\begin{matrix}X^{-1} \\Y\end{matrix}, Z^{-1})
\left\{ \begin{matrix} X^{-1} \\ Y \\ Z^{-1} \end{matrix} \right\} 
\\ &= \tau(A, A^h)\tau(\id_{l(A)}, A) \tau(A, A^h)^{-1} \, \sigma(A^{-1}, A^h)^{-1}
\tau(A, A^h)^{-1} \, \sigma(A^{-1}, A^h)^{-1} \\ &\qquad \qquad \times \sigma(A^{-1}, A^h) 
\sigma(\id_{b(A)^{-1}}, A^{-1}) \, A^{-1}\\
&= \tau(A, A^h)^{-1} \, \sigma(A^{-1}, A^h)^{-1} \, A^{-1} = \Ss(A),
\end{align*}
where the second sum is over all $X,Y,Z$ such that \begin{tabular}{p{0,4cm}|p{0,8cm}|p{0,4cm}}\quad & $X^{-1}$ & \quad
\\ \hline $X$ & $Y$ & $Z$ \\ \hline \quad  & $Z^{-1}$ & \quad  \end{tabular}
and $XYZ = A$; the third equality is by Lemma  \ref{l-atp3-bis};
and the fourth is clear.
 
\medbreak
(ii) is straightforward. \epf

\begin{prop}\label{hopf}
Let $\T$ be a vacant double groupoid and let $(\sigma, \tau)$ be a normalized
2-cocycle on $\T$ with values in $\ku^{\times}$.
Then $\ku{}_{\sigma}^{\tau} \T$ is a Hopf algebra if and only if $\T$ arises from a matched pair
of groups. \end{prop}

\pf Clearly, $\Delta(1) = 1 \otimes 1$ if and only if
$\sum_{A|B, \, AB = \iddv \,{t(AB)}} \tau (A, B) A \otimes B 
 = \sum_{x, y \in \Hc} \iddv \,{x} \otimes \iddv \,{y}$.
In particular, for any two arrows $x, y \in \Hc$, $\iddv \,{x} \vert \iddv \,{y}$.
This implies that $\Pc$ has only one point, \emph{i. e.} that $\T$ arises from a matched pair
of groups. The converse is well-known, so the proof is complete. \epf

\begin{prop}\label{involutory} 
Let $\T$ be a vacant double groupoid and let $(\sigma, \tau)$ be a normalized
2-cocycle on $\T$ with values in $\ku^{\times}$.
Then $\ku{}_{\sigma}^{\tau} \T$ is involutory.
If char $\ku = 0$, then $\ku{}_{\sigma}^{\tau} \T$ is semisimple and cosemisimple.
\end{prop}
\pf  Let $A\in \B$. We compute
$$
\Ss^2(A) = \tau(A, A^h)^{-1} \, \sigma(A^{-1}, A^h)^{-1}\, \Ss(A^{-1}) 
= \tau(A, A^h)^{-1} \, \sigma(A^{-1}, A^h)^{-1}\,
\tau(A^{-1}, A^v)^{-1} \, \sigma(A, A^v)^{-1}\, A.
$$
Now,
\begin{tabular}{p{0,5cm}|p{0,7cm}} $A$ & $A^h$
\\ \hline $A^v$ & $A^{-1}$ \end{tabular} implies that 
$$1 = \tau(A, A^h) \tau(A^v, A^{-1}) \sigma(A, A^v) \sigma(A^h, A^{-1}).$$
We conclude, using \eqref{cond-sigma-bis} and \eqref{cond-tau-bis},
that $\Ss^2(A) = A$. The second statement follows then from \cite[Corollary 6.5]{nik}.
\epf

The category of finite-dimensional representations of a weak Hopf algebra 
over $\ku$ admits the structure of a $\ku$-linear rigid monoidal category \cite{nik-v}.
Recall from \cite{ENO} the definition of a \emph{multifusion category}:
this is a semisimple $\ku$-linear rigid tensor category with finitely many isoclasses of simple objects 
and finite dimensional hom-spaces. 
Recall also that a multifusion category is called 
a \emph{fusion} category if in addition the unit object is \emph{simple}.

\begin{prop} Let $\T$ be a vacant double groupoid and let $(\sigma, \tau)$ be a normalized
2-cocycle on $\T$ with values in $\ku^{\times}$.

(i).  The unit object of the category $\Rep \ku{}_{\sigma}^{\tau} \T$ 
is simple if and only if  $\Vc \rightrightarrows \Pc$ is connected. 

\medbreak
(ii). If char $\ku = 0$, then  the category $\Rep \ku{}_{\sigma}^{\tau} \T$ of finite dimensional
$\ku{}_{\sigma}^{\tau} \T$-modules is a multifusion category. 
It is fusion if and only if  $\Vc \rightrightarrows \Pc$ is connected. 

\end{prop}

\pf We prove (i). We have already observed
that the target  subalgebra-- which is the unit object of $\Rep \ku{}_{\sigma}^{\tau} \T$ by
general reasons-- is the span of the elements ${}_{P}\uno$.
Let $\sim_V$ be the equivalence relation in $\Pc$ induced by $\Vc \rightrightarrows \Pc$.
We claim that the subspaces $\sum_{P \in X}\ku {}_{P}\uno$, for $X$ an equivalence class of $\sim_V$, are the simple subobjects of $\ku{}_{\sigma}^{\tau} \T_t$; whence the first equivalence. 
Indeed, for all $A\in \B$, we have 
$$A . {}_{P}\uno = \epsilon_t(A \, {}_{P}\uno) = \begin{cases} {}_{Q}\uno, \quad &\text{if $A = \id g$, for some $g \in \Vc$ such that $b(g) = P$, $t(g) = Q$,} \\
0, \quad &\text{otherwise}.
\end{cases}$$
This proves the claim. The second equivalence follows from Proposition \ref{hopf}.
Now (ii) follows from Proposition \ref{involutory}, the general theory of
weak Hopf algebras and (i).
\epf

\begin{prop} Let $\T$ be a vacant double groupoid and let $(\sigma, \tau)$ be a normalized
2-cocycle on $\T$ with values in $\ku^{\times}$.
Then $(\tau, \sigma)$ is a normalized 2-cocycle on the transpose
double groupoid $\T^t$ and the quantum groupoid
$\ku{}_{\tau}^{\sigma} \T^t$ is dual to $\ku{}^{\tau}_{\sigma} \T$.
\end{prop}
\pf The duality is given by the bilinear form $(B \vert C) =
\delta_{B, C^t}$, $B$, $C\in \B$. \epf

\subsection{The category $\Rep \ku \T$}

\

Let $\T$ be a finite vacant double groupoid, and let $\ku \T$ be the quantum groupoid associated to $\T$ as in Theorem \ref{bicross}. Consider the category $\C : = \Rep \ku \T$ of finite dimensional representations of $\ku\T$.  

Our aim in this subsection is to sketch a 
combinatorial description of the category  $\Rep \ku \T$ in groupoid-theoretical terms.  
We shall follow the lines in \cite[5.1]{nik-v}.

\bigbreak
Suppose that $s, e: \G \rightrightarrows \Pc$ is a groupoid. A \emph{$\ku$-linear $\G$-bundle}, 
or $\G$-bundle for short,
is a map $p: V \to \Pc$ together with an action of $\G$ on $p$, and such that 

(i) each fiber $V_b$ ($b \in \Pc$) is a vector space over $\ku$;

(ii) for all $g \in \G$ the map $g: V_{e(g)} \to V_{s(g)}$ is a linear isomorphism.

\medbreak
So one may think of a $\G$-bundle as a $\Pc$-graded vector space 
$V = \oplus_{b \in \Pc}V_b$ endowed with a linear $\G$-action $g: V_{e(g)} \to V_{s(g)}$, $g \in \G$.

\begin{obs}\label{equivcat} The category $\G$-bund of (finite dimensional) $\ku$-linear $\G$-bundles 
is equivalent to the category of (finite dimensional) representations of the groupoid 
algebra $\ku \G$. The equivalence is defined as follows: for a $\ku \G$-module $V$, 
we let the $\Pc$-grading on $V$ be given by $V_b = \id b . V$, for all $b \in \Pc$. \end{obs}

\bigbreak
Let now $\T$ be a vacant double groupoid. Let $\T$-bund be the category of 
$\ku$-linear bundles over the vertical groupoid $\B \rightrightarrows \Hc$.
Thus, the objects of  $\T$-bund are $\Hc$-graded vector spaces endowed with a left action of 
the vertical groupoid $\B \rightrightarrows \Hc$ by linear isomorphisms.

There is a structure of rigid monoidal category on $\T$-bund:

\begin{itemize}
\item \emph{Tensor product.} If $V$, $W$ are $\T$-bundles then $V \otimes W 
:= \oplus_{z\in \Hc} (V \otimes W)_z$, where $$(V \otimes W)_z 
= \sum_{xy = z}V_x \otimes_{\ku} W_y, \qquad z \in \Hc.$$
(Note that this difers from  $V \otimes_{\ku}W$ by the fact that we are not taking all summands $V_x \otimes_{\ku}W_y$ but only those for which $x$ and $y$ are composable.)
The action of $\B$ on $V \otimes W$ is given by $\Delta$.

\bigbreak
\item  \emph{Unit object.} This is the target subalgebra $\ku \T_t = \oplus_{P \in \Pc} \ku {}_P\uno$, with 
$\Hc$-grading defined by 
$$
(\ku \T_t)_x = \begin{cases} 0, \quad &\text{if $x$ is not an identity,} \\
\ku {}_{P}\uno, \quad &\text{if } x = \id P,
\end{cases} \quad \text{for all } x\in \Hc.
$$
and $\B$-action $A . {}_P\uno = \epsilon_t(A \, {}_P\uno)$.

\bigbreak
\item The \emph{dual} $V^*$ of an object $V = \oplus_{x \in \Hc}V_x \in \C$ has
$\Hc$-grading $(V^*)_x = (V_{x^{-1}})^*$, $x \in \Hc$;
$\B$-action $A : = (A^{-1})^* : (V^*)_{b(A)} \to (V^*)_{t(A)}$, for all $A \in \B$. \end{itemize}

\bigbreak
With  remark \ref{equivcat} in mind, we can describe the monoidal structure in $\Rep \ku \T$.  

\begin{prop} Let $\T$ be a vacant double groupoid. Assume that char $\ku =0$.
The category $\T$-bund  is a multifusion category over $\ku$ and it is monoidally 
equivalent to $\Rep \ku \T$.  \end{prop}

\pf The expressions for the tensor product, unit object and duals are  a translation of the formulas in \cite{nik-v} to the language of $\T$-bundles. For instance, the unit isomorphism $\ku \T_t \otimes V \to V$ is given as follows: for any $z \in \Hc$, we have  
$(\ku \T_t \otimes_{\ku}V)_z =  \ku {}_{l(z)}\uno \otimes_{\ku}V_z$; the isomorphism $\ku \T_t \otimes V \to V$ is determined by its homogeneous components $(\ku \T_t \otimes_{\ku}V)_z \to V_z$, given by the action of ${}_{l(z)}\uno$, which is the identity on $V_z$. 

The unit isomorphisms on the right  and the evaluation and coevaluation maps for the duals are constructed similarly. \epf

\subsection{A Kac exact sequence for matched pairs of groupoids}

\

We first recall the well-known definition of the groupoid cohomology via standard resolutions \cite{we, renault}.

Let $s, e:\G \rightrightarrows \Pc$ be a groupoid. In this subsection, we shall denote by
$\G^{(0)} : = \Pc$ the base of $\G$, $\G^{(1)} := \G$ and
$$\G^{(n)} =\{(x_1, \dots, x_{n}) \in \G^{n}:
x_1 \vert x_2 \vert \dots \vert x_{n-1}\vert x_{n} \}, \qquad n \ge 2.$$

Let $M$ be an abelian group (with trivial $\G$-action), and let 
$$C^n (\G, M) = \{f: \G^{(n)} \to M:  
f(x_1, \dots, x_{n}) = 0, 
\text{ if some }x_i\in \G^{(0)}\}.$$
The cohomology groups $H^n(\G, M)$ of $\G$ with coefficients in  $M$
 are the cohomology groups of the complex
\begin{multline}\label{complex}\begin{CD}
0 @>>> C^0 (\G, M) @>d^0>> C^1 (\G, M)
 @>d^1>> C^2 (\G, M)@>d^2>>
\dots \end{CD} \\
\begin{CD}
@>>> C^n (\G, M) @>d^n>> C^{n+1} (\G, M)
@>>> \dots
\end{CD}
\end{multline}
where 
\begin{equation}\label{cob} 
\begin{aligned}d^0 f(x) &= f(e(x)) - f(s(x)), \\
d^n f(x_1, \dots, x_{n+1}) &=  f(x_2, \dots, x_{n+1})
+ \sum_{1\le i \le n}(-1)^{i} f(x_1, \dots, x_ix_{i+1}, \dots, x_{n+1})
\\ & \qquad + (-1)^{n+1} f(x_1, \dots, x_{n}).
\end{aligned}
\end{equation}

\medbreak
Let now $\T$ be a  double groupoid: 
$\begin{matrix} \B &\rightrightarrows &\Hc
\\\downdownarrows &&\downdownarrows \\ \Vc &\rightrightarrows &\Pc \end{matrix}$. Let
\begin{align*}
\B^{(0,0)} &:= \Pc, \\
\B^{(0,s)} &:=  \left\{\left(x_1, \dots, x_s\right) \in  \Hc^{\, s}:
x_1\vert x_2 \dots \vert x_s \right\} = \Hc^{(s)}, \quad s> 0, \\
\B^{(r,0)} &:= \left\{\left(g_1, \dots, g_r\right) \in\Vc^{\, r}:
g_1\vert g_2 \dots \vert g_r\right\} = \Vc^{(r)}, \quad r> 0,\\
\B^{(r,s)} &:= \left\{\left(\begin{tabular}{p{0,8cm} p{0,8cm} p{0,8cm} p{0,8cm}}
$A_{11}$ & $ A_{12}$ & \dots & $A_{1s}$ \\
$A_{21}$ & $ A_{22}$ & \dots & $A_{2s}$ \\
\dots & \dots & \dots & \dots \\ 
$A_{r1}$ & $ A_{r2}$ & \dots & $A_{rs}$  \end{tabular}\right) \in 
\B^{r\times s}: \quad
\begin{tabular}{p{0,8cm}|p{0,8cm}|p{0,8cm}|p{0,8cm}}
$A_{11}$ & $ A_{12}$ & \dots & $A_{1s}$ \\ \hline
$A_{21}$ & $ A_{22}$ & \dots & $A_{2s}$ \\ \hline
\dots & \dots & \dots & \dots \\ \hline
$A_{r1}$ & $ A_{r2}$ & \dots & $A_{rs}$  \end{tabular} \right\}, 
\quad r,s> 0.
\end{align*}

\bigbreak
Let $M$ be an abelian group and let $D^{r, s} = D^{r, s}(\T, M)$, $r, s \geq 0$, be defined by
\begin{equation*}
D^{r, s}  := \left\{f: \B^{(r,s)} \to M:  
f\left(\begin{tabular}{p{0,8cm}  p{0,8cm} p{0,8cm}}
$A_{11}$  & \dots & $A_{1s}$ \\
$A_{21}$  & \dots & $A_{2s}$ \\
\dots  & \dots & \dots \\ 
$A_{r1}$  & \dots & $A_{rs}$  \end{tabular} \right) 
= 0,  \text{ if}  \begin{cases} \text{either  } r > 1, s>0  \text{ and } A_{ij} \in \Vc, \\
\text{or  } r > 1, s = 0  \text{ and } A_{i0} \in \Pc,  \\
 \text{or  } r>0, s > 1  \text{ and } A_{ij} \in \Hc,\\
\text{or  } r =0, s > 1  \text{ and } A_{0j} \in \Pc.  \end{cases}
\right\}.
\end{equation*}

Let $d_H = d_H^{r,s}: D^{r, s} \to D^{r, s+1}$, 
$d_V = d_V^{r,s}: D^{r, s} \to D^{r+1, s}$ 
be, respectively, the horizontal and vertical coboundary maps defined as follows: 

\begin{itemize}
\item If  $r =0$, $d_H$ is as in \eqref{cob};

\item if  $s =0$, $d_V$ is as in \eqref{cob};

\item if  $r =0$, $s > 0$,
$$d_V^{0,s} f(A_{11}, \dots, A_{1s}) = 
f(b(A_{11}), \dots, b(A_{1s})) - f(t(A_{11}), \dots, t(A_{1s}));$$
\item if  $r >0$, $s = 0$,
$$d_H^{r,0} f\begin{pmatrix}
A_{11}   \\
\dots   \\ 
A_{r1}  \end{pmatrix} = 
f(r(A_{11}), \dots, r(A_{r1})) - f(l(A_{11}), \dots, l(A_{r1}));$$
\item if $r > 0$ and $s > 0$,
\begin{align*} 
d_V^{r,s} f\begin{pmatrix}
A_{11}  & \dots & A_{1s} \\
\dots  & \dots & \dots \\ 
A_{r1}  & \dots & A_{rs} \\
A_{r+1,1}  & \dots & A_{r+1,s} 
\end{pmatrix} &=  
f\begin{pmatrix}
A_{21}  & \dots & A_{2s} \\
\dots  & \dots & \dots \\ 
A_{r1}  & \dots & A_{rs} \\
A_{r+1,1}  & \dots & A_{r+1,s} 
\end{pmatrix}
\\ & \qquad + \sum_{1\le i \le r}(-1)^{i} f\begin{pmatrix}
A_{11}  & \dots & A_{1s} \\
\dots  & \dots & \dots \\ 
\left\{\begin{matrix}A_{i1} \\A_{i+1,1}\end{matrix}\right\}  
& \dots & \left\{\begin{matrix}A_{is} \\A_{i+1,s}\end{matrix}\right\} \\
\dots  & \dots & \dots \\ 
A_{r+1,1}  & \dots & A_{r+1,s} 
\end{pmatrix}
\\ & \qquad + (-1)^{r+1} f\begin{pmatrix}
A_{11}  & \dots & A_{1s} \\
\dots  & \dots & \dots \\ 
A_{r1}  & \dots & A_{rs} \end{pmatrix};
\end{align*} 
\begin{align*} 
d_H^{r,s} f\begin{pmatrix}
A_{11}  & \dots  & A_{1,s+1} \\
\dots  & \dots & \dots \\ 
A_{r1}  & \dots & A_{r,s+1} \end{pmatrix} &=  
f\begin{pmatrix}
A_{12}  & \dots & A_{1,s+1} \\
\dots  & \dots & \dots \\ 
A_{r2}  & \dots & A_{r,s+1}  
\end{pmatrix}
\\ & \qquad + \sum_{1\le j \le s}(-1)^{j} 
f\begin{pmatrix}
A_{11}  & \dots & \left\{A_{1,j} A_{1,j+1} \right\}& \dots & A_{1,s+1} \\
\dots  & \dots& \dots& \dots & \dots \\ 
A_{r,1}  & \dots& \left\{A_{r,j}A_{1,j+1} \right\}& \dots & A_{r,s+1} 
\end{pmatrix}
\\ & \qquad + (-1)^{s+1} f\begin{pmatrix}
A_{11}  & \dots & A_{1s} \\
\dots  & \dots & \dots \\ 
A_{r1}  & \dots & A_{rs} \end{pmatrix}.
\end{align*}
\end{itemize}

A straightforward computation shows that the following diagram
commutes:
$$
\begin{CD}
D^{r+1, s} @>{d_H}>> D^{r+1, s+1}
\\
@AA{d_V}A @AA{d_V}A \\
D^{r, s} @>{d_H}>> D^{r, s+1}
\end{CD} 
$$

Thus, there is a double cochain complex 
$$ D^{\cdot\cdot} = \qquad
\begin{CD}
\vdots \\
@AAA \\
D^{2, 0} @>{d_H}>> \hspace{18pt}\raisebox{1.0ex}{\vdots}\cdots \\
@AA{d_V}A @AA{-d_V}A \\
D^{1, 0} @>{d_H}>> D^{1, 1} @>{d_H}>>
\hspace{18pt}\raisebox{1.0ex}{\vdots}\cdots \\
@AA{d_V}A @AA{-d_V}A @AAA \vspace{9pt}\\
D^{0, 0} @>{d_H}>> D^{0, 1} @>{d_H}>>
D^{0, 2} @>>> \cdots\ \ ,
\end{CD} 
$$
with the usual "sign trick": the vertical differential
is $(-1)^{s} d_V^{r,s}$. We then remove the edges
of this double complex setting $A^{r, s}(\T, M) = A^{r,s} := D^{r+1,s+1}$,
$r,s \ge 0$; and denote by $E^{\cdot\cdot}(\T, M) = E^{\cdot\cdot}$ the double complex 
consisting only
of the edges of $D^{\cdot\cdot}$. Compare with  
\cite[pp. 173 ff.]{ma-newdir}.

\bigbreak
We are now ready to state a result inspired in the celebrated
Kac exact sequence \cite[(3.14)]{k}. Suppose that
$\T$ is a \emph{vacant} double groupoid
and let  $\D = \Vc \bowtie \Hc$ be corresponding diagonal groupoid 
as in Proposition \ref{equiv-matchedpair}.

\begin{prop}\label{kes} 
There is an exact sequence 
\begin{equation}\label{kes-formulagral}
\begin{aligned}
0 &\to H^1(\D, M)\to H^1(\Hc, M) \oplus H^1(\Vc, M) \to 
H^0(\Tot A(\T, M)^{\cdot\cdot}, M) \\
&\to H^2(\D, M)\to H^2(\Hc, M) \oplus H^2(\Vc, M) \to 
H^1(\Tot A(\T, M)^{\cdot\cdot}, M) \\
&\to H^3(\D, M)\to H^3(\Hc, M) \oplus H^3(\Vc, M).\end{aligned}
\end{equation}
\end{prop}

\pf The short exact sequence of double complexes 
$\begin{CD} 0@>>> A^{\cdot\cdot} @>>> D^{\cdot\cdot}
@>>> E^{\cdot\cdot}@>>> 0\end{CD}$ 
(with $A^{\cdot\cdot}$ "shifted") induces a long exact sequence in cohomology.
It is clear that 
\begin{equation*}
H^n(\Tot E^{\cdot\cdot} (\T, M)) =  H^n(\Hc, M) 
\oplus H^n(\Vc, M), \quad n> 0.
\end{equation*}
We claim that
\begin{equation}\label{kes-aux}
H^n(\Tot D^{\cdot\cdot} (\T, M))) =  H^n(\D, M), \quad n> 0.
\end{equation}
Indeed, $H^{\cdot}(\D, M)$ are the cohomology groups of
a complex $\Hom_{\ku \D} (F^{\cdot}, M)$, where $F^{\cdot}$
is some free resolution of the trivial $\D$-module.
Now, arguing as in \cite[Lemma 1.7]{ma-newdir}, 
we see that $H^{\cdot}(\Tot D^{\cdot\cdot} (\T, M)))$
are also the cohomology groups of
a complex $\Hom_{\ku \D} (G^{\cdot}, M)$, where $G^{\cdot}$
is another free resolution of the trivial $\D$-module;
this implies \eqref{kes-aux}. \epf

If $M = \ku^{\times}$, it is natural to denote
\begin{align}\Aut(\ku \T) &= H^0(\Tot A^{\cdot\cdot}(\T,\ku^{\times})),\\
\Opext(\ku^{\Vc}, \ku\Hc) &= H^1(\Tot A^{\cdot\cdot} (\T,\ku^{\times})),
\end{align}
by Theorem \ref{concociclos} and 
in view of an extension theory of quantum groupoids yet to be
explored. Then \eqref{kes-formulagral} has in this case the familiar
expression
\begin{equation}\label{kes-formula}
\begin{aligned}
0 &\to H^1(\D, \ku^{\times})\to H^1(\Hc, \ku^{\times}) \oplus H^1(\Vc, \ku^{\times}) \to \Aut(\ku \T) \\
&\to H^2(\D, \ku^{\times})\to H^2(\Hc, \ku^{\times}) \oplus H^2(\Vc, \ku^{\times}) \to \Opext(\ku \T) \\
&\to H^3(\D, \ku^{\times})\to H^3(\Hc, \ku^{\times}) \oplus H^3(\Vc, \ku^{\times}).
\end{aligned}
\end{equation}

\subsection{Conclusion}

\

We have introduced  families of quantum groupoids and \emph{a fortiori}
of tensor categories. To be sure that these tensor categories
are really new, we have to explicitly compute first the 
$\Opext(\ku \T)$ groups, and second to analyze when the corresponding
quantum groupoids give rise to equivalent tensor categories.
We shall address both questions in subsequent work.

\end{document}